\documentclass{article}
\usepackage{amsmath}
\usepackage{amssymb}
\usepackage{amsthm}
\usepackage{stmaryrd}
\usepackage{wasysym}
\usepackage{tikz-cd}
\newtheorem{theorem}{Theorem}[section]
\newtheorem{definition}[theorem]{Definition}
\newtheorem{proposition}[theorem]{Proposition}
\newtheorem{lemma}[theorem]{Lemma}
\newtheorem{corollary}[theorem]{Corollary}
\newtheorem{remark}[theorem]{Remark}
\usepackage{indentfirst}
\usepackage[english]{babel}
\usepackage[letterpaper,top=2cm,bottom=2cm,left=3cm,right=3cm,marginparwidth=1.75cm]{geometry}
\usepackage{amsmath}
\usepackage{graphicx}
\usepackage[colorlinks=true, allcolors=blue]{hyperref}
\title{Toward the profinite rigidity of hyperbolic 3-manifolds}
\author{Tianwei Liu}
\begin{document}
\date{}
\maketitle
\begin{abstract}
This note surveys recent progress toward the profinite rigidity of orientable finite-volume hyperbolic 3-manifolds. Beginning in a brief review of some basic settings of profinite completion and rigidity of general groups, we state the complete results in profinite detection of 3-manifold decompositions and profinite rigidity or nonrigidity of other types of 3-manifolds, which reduce the problem to the hyperbolic case. Then we give some evidence to a positive answer and conclude some existing ideas toward the remained hyperbolic case. We also summarize the methods and techniques used in the process.
\end{abstract}
\tableofcontents
\section{Introduction}
In this paper, the author records recent progress toward the profinite rigidity of orientable finite-volume hyperbolic 3-manifolds. This paper also summarizes the techniques and ideas developed in the process of research. \par
The organization of the paper is as below. Section 2 is aimed for a brief introduction to settings of profinite completion and rigidity for a general group. Section 3 presents the process to reduce the problem to the hyperbolic case, including the profinite detection of prime decomposition and JSJ decomposition of 3-manifolds due to Wilton and Zalesskii in \cite{WltZ19} and the satisfactory solution to the other seven geometric cases other than the hyperbolic case. In Section 4, we concentrate on the hyperbolic case and present some evidence toward a positive answer, including some examples of profinite rigid hyperbolic 3-manifolds (see \cite{BRW17} and \cite{Che24} for example), the profinite almost rigidity of 3-manifolds (see Table 1 and Table 2 of \cite{Xu25} for a summary and \cite{Liu23} for the proof of the hyperbolic case) and the more recent work by Agol, Cheetham-West and Minsky \cite{ACM24}. Sketch of Proofs and justification of concepts are delayed to Section 5 as a summary of methods. 
\section*{Acknowledgement}
The author claims no originality to the material and are thankful to the authors of related papers for reference.  
\section{Preliminaries for profinite settings}
To begin with we briefly review some basic settings in profinite completion and rigidity of general groups. For a thorough and readable introduction the reader is referred to \cite{Rei18}. When studying a certain group, it is natural to ask that to what extent do its finite quotients determine the group itself, at least among some specified family of groups. To be precise, for any finitely generated group $\Gamma$, we define
\begin{displaymath}
    \widehat{\Gamma}=\varprojlim_{N\lhd\Gamma\ of\ finite\ index}\Gamma/N
\end{displaymath}
to be the \textit{profinite completion} of $\Gamma$. Namely, we take the directed set $\mathcal{N}$ which consists of normal subgroups of $\Gamma$ of finite index. For any $N$ and $N'\in\mathcal{N}$, set $N>N'$ if $N\subset N'$. The homomorphism $\Gamma
/N\rightarrow\Gamma/N'$ is given by the natural quotient. Note that finitely generation of $\Gamma$ ensures that the number of normal subgroups of bounded indices is finite. Therefore, we can construct a natural cofinal sequence $\{\Gamma_n\}_{n\in\mathbb{N}_+}$ of characteristic subgroups of the colimit where
\begin{displaymath}
    \Gamma_n=\bigcap_{N\lhd\Gamma\ with\ [\Gamma:N]\leqslant n }N,
\end{displaymath}
which makes things more convenient in many occasions. There is a natural homomorphism $\iota:\Gamma\rightarrow\widehat{\Gamma}$ which maps $g\in\Gamma$ to $(gN)_{N\lhd\Gamma\ of\ finite\ index}$. It is easily qualified that $\iota$ is injective if and only if $\Gamma$ is residually finite. Unless otherwise stated, we will assume finitely generation and residually finiteness for all groups discussed since the fundamental groups of compact 3-manifolds are always finitely generated and residually finite (the former is well known while for the latter we refer to \cite{Hem87} and Theorem 3.3 of \cite{Thu82}). \par
Now, it is not difficult to see that for any finite generated residually finite groups $\Gamma_1$, $\Gamma_2$, $\widehat{\Gamma_1}\cong\widehat{\Gamma_2}$ if and only if for any finite quotient of $\Gamma_1$(respectively, $\Gamma_2$), there exists a finite quotient of $\Gamma_2$(respectively, $\Gamma_1$) which is isomorphic to it, that is to say, the profinite completion of a group organizes total information of all of its finite quotients.\par
We call a finitely generated residually finite group $\Gamma$ \textit{profinitely rigid} if for any finitely generated residually finite group $\Gamma'$ such that $\widehat{\Gamma}\cong\widehat{\Gamma'}$, we always have $\Gamma\cong\Gamma'$. In other word, a group is profinitely rigid if its profinite completion (or the set of its finite quotients) determines itself uniquely. Let $\mathcal{C}$ be a family of groups, such as the set of fundamental groups of compact 3-manifolds. For any $\Gamma\in\mathcal{C}$, we say $\Gamma$ is \textit{profinitely rigid in $\mathcal{C}$} if for any $\Gamma'\in\mathcal{C}$, $\widehat{\Gamma}\cong\widehat{\Gamma'}$ implies $\Gamma\cong\Gamma'$.\par
We give $\widehat{\Gamma}$ a topological space structure by giving each finite quotient $\Gamma/N$ a discrete topology and equipping $\Gamma$ with the inverse limit topology. More generally, a \textit{profinite space} (respectively, a \textit{profinite group}) is an inverse limit of a family of finite spaces (respectively, finite groups) each equipped with the discrete topology over a directed partially ordered set. It is shown in \cite{RibZ10}, Theorem 1.1.12 that $X$ is a profinite space if and only if $X$ is compact, Hausdorff and totally disconnected, while it is shown in \cite{RibZ10}, Theorem 2.1.3 that $G$ is a profinite group if and only if $G$ is a topology group whose underlying space is a profinite space. For more about profinite groups, see \cite{RibZ10}. Therefore the profinite completion of a group is a profinite topology group. A deep theorem in \cite{NS07} shows that any subgroup of finite index in a finitely generated profinite group is open. In general, we can also define the profinite topology on any finitely generated residually finite group $\Gamma$ by assigning the open sets to be its finite index normal subgroups and their cosets. The profinite completion of $\Gamma$ is just the completion with respect to its profinite topology. \par
For subgroups of $\Gamma$ and $\widehat{\Gamma}$, we have the following.
\begin{proposition}[\cite{Rei18}, Proposition 2.1]
    If $\Gamma$ is a finitely generated residually finite group, then there is a one-to-one correspondence between the set $\mathcal{X}$ of subgroups of $\Gamma$ that are open in the profinite topology on $\Gamma$, and the set $\mathcal{Y}$ of all finite index subgroups of $\widehat{\Gamma}$.\par
    Identifying $\Gamma$ with its image in the completion, this correspondence is given by
    \begin{itemize}
        \item For $H\in\mathcal{X}$, $H\mapsto\overline{H}$.
        \item For $Y\in\mathcal{Y}$, $Y\mapsto Y\cap\Gamma$.
    \end{itemize}
    If $H,K\in\mathcal{X}$ and $K<H$ then $[H:K]=[\overline{H}:\overline{K}]$. Moreover, $K\triangleleft H$ if and only if $\overline{K}\triangleleft\overline{H}$, and $\overline{H}/\overline{K}\cong H/K$.
\end{proposition}
Note that for any finitely generated residually finite group $\Gamma$ and its subgroup $H<\Gamma$, the profinite topology on $\Gamma$ determines some pro-topology on $H$ and therefore some completion of $H$. Since $\Gamma$ is residually finite, $H$ injects into $\widehat{\Gamma}$ and determines a subgroup $\overline{H}\subseteq\widehat{\Gamma}$. Hence there is a natural epimorphism $\widehat{H}\rightarrow\overline{H}$. This need not be injective. For this to be injective (that is, the full profinite topology is induced on $H$), it is easy to see that the following needs to hold: \textit{For every subgroup $H_1$ of finite index in $H$, there exists a finite index subgroup $\Gamma_1<\Gamma$ such that $\Gamma_1\cap H<H_1$}.\par
With the same assumption, we call $\Gamma$ \textit{H-separable} if for every $g\in G\backslash H$, there is a subgroup $K$ of finite index of $\Gamma$ such that $H\subseteq K$ but $g\notin K$, or equivalently, the intersection of all finite index subgroups in $\Gamma$ containing $H$ is precisely $H$. We call $\Gamma$ \textit{LERF} (or subgroup separable) if it is $H$-separable for every finitely generated subgroup $H$, or equivalently, if every finitely generated subgroup is a closed subset in the profinite topology. If $\Gamma$ is LERF, then for any finitely generated subgroup $H$, the profinite topology on $\Gamma$ induces the full profinite topology on $H$, or equivalently, the epimorphism $\widehat{H}\rightarrow\overline{H}$ is an isomorphism (see \cite{Rei18}, Lemma 2.8 and Corollary 2.9). 
\section{Reduction to hyperbolic manifolds}
In this section, we clarify what the problem of profinite rigidity of 3-manifolds is and how it reduces to the hyperbolic case, as well as the complete results in other types of 3-manifolds.\par
A 3-manifold $M$ is called \textit{profinitely rigid} if $\pi_1M$ is profinitely rigid among the family of fundamental groups of 3-manifolds, that is, for any 3-manifold $N$, $\widehat{\pi_1M}\cong\widehat{\pi_1N}$ implies $\pi_1M\cong \pi_1N$ (note that in general, spherical 3-manifolds and bounded 3-manifolds are not uniquely determined by their fundamental groups, so we do not require $M\cong N$). Otherwise we say that it is \textit{profinitely nonrigid}. In the discussion below, we refer to a 3-manifold for a connected, compact, orientable one with empty or toroidal boundary (see, for example, \cite{AscFW15}, Lemma 1.4.3 and the discussion there). 
\subsection{Profinite detection of 3-manifold decompositions}
\textbf{3.1.1. Profinite detection of prime decompositions.} For any 3-manifold $M$, there exists a unique \textit{prime decomposition} representing it to connected sum of $S^1\times S^2$'s and irreducible pieces, that is,
\begin{displaymath}
    M=m(S^1\times S^2)\#N_1\#\cdots\#N_n
\end{displaymath}
where $m,n\in\mathbb{N}$ and the $N_i$'s are irreducible. Wilton and Zalesskii have proved in \cite{WltZ19} that the profinite completion of its fundamental group determines any 3-manifold's prime decomposition. 
\begin{theorem}[\cite{WltZ19}, Theorem A]
    Let $M,N$ be closed orientable 3-manifolds with prime decompositions $M\cong M_1\#\cdots\#M_m\#r(S^1\times S_2)$ and $N\cong N_1\#\cdots\#N_n\#s(S^1\times S_2)$. If $\widehat{\pi_1M}\cong\widehat{\pi_1N}$, then $m=n,\ r=s$, and up to re-ordering the image of $\widehat{\pi_1M_i}$ is conjugate to $\pi_1N_i$ for each $i$.
\end{theorem}
The theorem above says that to solve the problem of profinite rigidity of 3-manifolds, it suffices to concentrate on irreducible 3-manifolds.\par
\textbf{3.1.2. Profinite detection of JSJ decompositions. }Next, for any irreducible 3-manifold, there exists a unique JSJ decomposition which divides it into Seifert fibered pieces and atoroidal pieces by a (possibly empty) collection of disjoint embedded incompressible tori. Namely, we have:
\begin{theorem}[JSJ Decomposition Theorem; \cite{AscFW15}, Theorem 1.6.1]
    Let $M$ be a compact, orientable, irreducible 3-manifold with empty or toroidal boundary. There exists a (possibly empty) collection of disjointly embedded incompressible tori $T_1,\dots ,T_m$ such that each component of N cut along $T_1,\dots,T_m$ is atoroidal or Seifert fibered. Any such collection of tori with a minimal number of components is unique up to isotopy.
\end{theorem}
Wilton and Zalesskii also proved in \cite{WltZ19} that the profinite completion of an irreducible 3-manifold determines its JSJ decomposition. However, we need more definitions and constructions to state it.\par
At this point we only give an intuitive understanding. A precise statement and a more detailed treatment is delayed to Subsection 5.1, see Theorem 5.6, invoking the Bass-Serre theory in geometric group theory. The reader is also referred to Theorem 4.3 of \cite{WltZ19} for the original material and Section 3 of \cite{Xu25}. We construct a graph $\mathcal{G}=(\mathcal{V},\mathcal{E})$ corresponding to the given JSJ decompositions of an irreducible 3-manifold $M$. The vertex set $\mathcal{V}=\{M_1,\cdots,M_m\}$ consists of the JSJ pieces of $M$ and the edge set $\mathcal{E}=\{T_1,\cdots,T_n\}$ consists of the JSJ tori of $M$. An edge $T_i$ connect two vertices $M_j,M_k$ if and only if $T_i$ is a public boundary component for both. Then Theorem 3.2 implies that the corresponding graphs of $M$ and $N$ are isomorphic and the fundamental groups of corresponding JSJ pieces are profinitely isomorphic, that is, $\widehat{\pi_1M_i}\cong\widehat{\pi_1N_i}$. Furthermore, the pattern of connection of the corresponding vertices and edges, $\widehat{\pi_1T_{M,i}}\rightarrow\widehat{\pi_1M_j}$ and $\widehat{\pi_1T_{N,i}}\rightarrow\widehat{\pi_1N_j}$ are natural under the graph isomorphism as well. All of the details are located in Section 5.1. \par
\begin{theorem}[\cite{WltZ19}, Theorem B]
    Let $M$ and $N$ be close, orientable, irreducible 3-manifolds, and suppose $\widehat{\pi_1M}\cong\widehat{\pi_1N}$. Then the underlying graphs of the JSJ decompositions of $\pi_1M,\pi_1N$ are isomorphic, and corresponding vertex groups have isomorphic profinite completions.
\end{theorem}
We remark that for the simpler graph manifold case the same authors proved this result much earlier in \cite{WltZ10}.\par
Theorem 3.1 and 3.3 altogether imply that in the research of the profinite rigidity of 3-manifolds, it suffices to consider the case when they are irreducible and are either Seifert fibered or atoroidal.\par
For a more detailed introduction of decomposition theorems of 3-manifolds, see Section 1 of \cite{AscFW15}. We refer to \cite{AscFW15} for a general introduction of recent progress of the study of 3-manifolds and their fundamental groups as well.
\begin{remark}[Distinguishing geometry]
    \rm{Before we enter into the next subsection, it is worthwhile to note that Wilton and Zalesskii have provided in another paper \cite{WltZ17} a profinite group-theoretical characterization of hyperbolic 3-manifolds and Seifert fibered 3-manifolds, which can be seen as some profinite analogue of the classical Hyperbolization Theorem and Seifert Conjecture (now theorem). We will revisit them in Section 5 as examples of profinite extension of classical results in 3-manifold topology. To be precise, we record the next results.}
\end{remark}
\begin{theorem}[\cite{WltZ17}, Theorem A]
    Let $M$ be a closed, orientable, aspherical 3-manifold. Then $M$ is hyperbolic if and only if the profinite completion $\widehat{\pi_1M}$ does not contain a subgroup isomorphic to $\widehat{\mathbb{Z}}^2$.
\end{theorem}
Note that Theorem 3.5 implies that any 3-manifold profinitely isomorphic to a hyperbolic manifold also yields a hyperbolic geometry. 
\begin{theorem}[\cite{WltZ17}, Theorem B]
    Let $M$ be a closed, orientable, aspherical 3-manifold. Then $M$ is Seifert fibered if and only if the profinite completion $\widehat{\pi_1M}$ has a non-trivial procyclic normal subgroup.
\end{theorem}
\subsection{Results in non-hyperbolic cases}
We begin with the Seifert fibered case. The theorem below shows the possible geometric structures that a Seifert fibered 3-manifold may have. Here $K^2\widetilde{\times}I$ denotes the unique orientable interval bundle over the Kleinian bottle $K^2$. 
\begin{proposition}[Classification of Seifert fibered 3-manifolds; \cite{AscFW15}, Theorem 1.8.1]
    Let $N$ be a compact, orientable 3-manifold with empty or toroidal boundary, which is not homeomorphic to $S^1\times D^2$, $T^2\times I$, or $K^2\widetilde{\times}I$. Then $N$ is Seifert fibered if and only if $N$ admits a geometric structure modeled on one of the following geometries: $\mathbb{S}^3$, $\mathbb{E}^3$, $S^2\times\mathbb{E}^1$, $\mathbb{H}^2\times\mathbb{E}^1$, $\widetilde{SL_2(\mathbb{R})}$, or $Nil$. 
\end{proposition}
We begin with the easiest spherical geometric case. Since the profinite completion of a finite group is itself, it is trivial that finite groups are profinitely rigid. Then the following Elliptization Theorem shows that spherical 3-manifolds are profinitely rigid. 
\begin{theorem}[Elliptization Theorem; \cite{AscFW15}, Theorem 1.7.3]
    3-manifolds with finite fundamental groups is spherical. 
\end{theorem}
For other Seifert fibered geometries, Hempel \cite{Hem14} shows that closed 3-manifolds with $\mathbb{H}^2\times\mathbb{E}^1$ geometry are in general profinitely nonrigid. Wilkes \cite{Wilk17} shows that closed 3-manifolds with $\mathbb{E}^3$, $Nil$, $S^2\times\mathbb{E}^1$, or $\widetilde{SL_2(\mathbb{R})}$ geometry are profinitely rigid. Finally, Corollary 8.3 of \cite{Xu25} shows that bounded Seifert fibered 3-manifolds are profinitely rigid.\par
The Hyperbolic Theorem below equates the atoroidal case to the hyperbolic case.
\begin{theorem}[Hyperbolization Theorem; \cite{AscFW15}, Theorem 1.7.5]
    Let $N$ be a compact, orientable, irreducible 3-manifold with empty or toroidal boundary. Suppose that N is atoroidal and not homeomorphic to $S^1\times D^2$, $T^2\times I$, or $K^2\widetilde{\times}I$. If $\pi_1N$ is infinite, then $N$ is hyperbolic. 
\end{theorem}
Therefore, after clearing finitely many trivial cases ($S^1\times D^2$, $K^2\widetilde{\times}I$, etc.), it suffices to research whether hyperbolic 3-manifolds are profinitely rigid, which is the task of Section 4. However, before we turn to Section 4, it is worthwhile to sketch some discussion in three special cases where the manifold has nontrivial JSJ decomposition in the next subsection, namely, $Sol$ manifolds, graph manifolds and mixed manifolds.\par
\subsection{Results in 3-manifolds with nontrivial JSJ decompositions}
Despite of (but making use of, certainly) the strength of decomposition methods, there are also some results for 3-manifolds with nontrivial JSJ decompositions. Namely, we discuss $sol$ geometric manifolds, graph manifolds and mixed manifolds. Note that a complete solution to the last case depends on the solution to the hyperbolic case.\par
\textbf{3.3.1. $Sol$ manifolds. }To be precise, we begin with a characterization of $Sol$ manifolds, which is a combination of Proposition 1.9.2 and Proposition 1.10.1 of \cite{AscFW15}.
\begin{proposition}
    Let $M$ be a 3-manifold other than $S^1\times D^2$, $T^2\times I$ and $K^2\widetilde{\times}I$. Then $M$ is of a $Sol$ geometry if and only if one of the followings holds. \par
    (1) $M$ is a torus bundle with an Anosov monodromy; \par
    (2) $M$ is a twisted double of $K^2\widetilde{\times}I$. 
\end{proposition}
Therefore a general $Sol$ manifold is an Anosov torus bundle, which is shown in \cite{Fun13} and \cite{Ste72} to be profinite nonrigid in general.\par
\textbf{3.3.2. Graph manifolds. }Now we turn to graph manifolds. An irreducible 3-manifold $M$ is called a \textit{graph manifold} if all of its JSJ pieces are Seifert fibered. Graph manifolds have a close connection to geometric group theory, in particular the research of their underlying JSJ graphs as in Section 3.1. Wilkes gave a complete solution to the profinite properties of graph manifolds. Namely, he proved in \cite{Wilk18} that a graph manifold is profinitely rigid if the underlying JSJ graph is nonbipartite. Otherwise, it is in general profinitely nonrigid. Recall that a graph $G=(V,E)$ is called \textit{bipartite} if there exists a division $V=V_1\bigsqcup V_2$ of its vertex set, such that for any $(v_i,v_j)$ in $E$, $v_i\in V_1,v_j\in V_2$ or $v_j\in V_1,v_i\in V_2$. The precise statement is as below.
\begin{theorem}[\cite{Wilk18}, Theorem 10.9]
    Let $M$ and $N$ be graph manifolds with JSJ decomposition graphs $X$ and $Y$ respectively.\par
    (1) If $X\cong Y$ is not bipartite, then $M$ is profinitely rigid;\par
    (2) If $X\cong Y$ is bipartite on two sets $R$ and $B$, then $\pi_1M$ and $\pi_1N$ have isomorphic profinite completions if and only if for some choices of generators of fibered subgroups, there is a graph isomorphism $\phi:\ X\rightarrow Y$ and some $\kappa\in\widehat{\mathbb{Z}}^{\times}$ such that:\par
    (a) For each edge $e$ of $X$, $\gamma(\phi(e))=\pm_e\gamma(e)$, where the sign is positive if both end vertices of $e$ have orientable base.\par
    (b) The total slope of every vertex space of $M$ or $N$ vanishes.\par
    (c1) If $d_0(e)=r\in R$, $\delta(\phi(e))=\pm_e\kappa\delta(e)$ modulo $\gamma(e)$, and $(M_r,N_{\phi(r)})$ is a Hempel pair of scale factor $\kappa$.
    (c2) If $d_0(e)=b\in B$, $\delta(\phi(e))=\pm_e\kappa^{-1}\delta(e)$ modulo $\gamma(e)$, and $(M_b,N_{\phi(b)})$ is a Hempel pair of scale factor $\kappa^{-1}$.
\end{theorem}
Clarifying Theorem 3.7 more precisely is beyond the scope of this paper. See the original paper \cite{Wilk18} for a detailed introduction. See also Section 5.1 for some background settings.\par
\textbf{3.3.3. Mixed manifolds. }We refer to a \textit{mixed manifold} as an irreducible 3-manifold which has JSJ decomposition with at least one hyperbolic piece. Since the profinite rigidity of hyperbolic 3-manifolds is still unclear, we certainly cannot expect a satisfactory result in mixed manifolds. However, benefiting from the solution to the Seifert fibered case, \cite{Xu25}, Corollary B has proved that the \textit{Seifert part} of a closed mixed manifold is profinitely rigid. In other word, whether a closed mixed manifold is profinitely rigid depends only on its hyperbolic pieces. However, Xu also proved that a bounded mixed manifold is in general profinitely nonrigid. In the argument the next theorem is important:
\begin{theorem}[\cite{Xu25}, Theorem A]
    Suppose $M$ and $N$ are mixed 3-manifolds, with empty or toroidal boundary. If $\widehat{\pi_1M}\cong\widehat{\pi_1N}$ through an isomorphism respecting the peripheral structure, then the Seifert parts of $M$ and $N$ are homeomorphic.
\end{theorem}
The result in closed mixed manifolds then immediately follows since by definition profinite isomorphism between closed manifolds respects the trivial peripheral structure. We will turn to the peripheral structure in Section 5.  
\section{Evidence to a positive answer}
In this section, we provide some evidence to a positive answer to the hyperbolic case, including some verified profinitely rigid families of hyperbolic manifolds and a complete solution to the problem of almost rigidity. 
\subsection{Some examples of rigidity}
One-punctured torus bundles and four-punctured sphere bundles are verified to be profinitely rigid (\cite{BRW17} and \cite{Che24}). These are all fibered hyperbolic 3-manifolds over the circle, which is simpler but essential, see Section 5.3 for more information. The simplest surface of negative Euler characteristic is $S^3\backslash3pts$, the pair of pants, on which fibered bundles are too simple to consider because of the trivial structure of its mapping class group $Mod(S^3\backslash3pts)=\Sigma_3$. A little more complex are $T^2\backslash pt$ and $S^3\backslash4pts$. However, they still have explicitly-presenting mapping class groups $Mod_\pm(T^2\backslash pt)=Out(F_2)=GL_2(\mathbb{Z})$ and $Mod(S^3\backslash4pts)=(\mathbb{Z}/2\mathbb{Z}\times\mathbb{Z}/2\mathbb{Z})\rtimes PSL_2(\mathbb{Z})$. For an introduction on mapping class groups we refer to \cite{FarM12} and see Chapter 2 for the presentations above. Note that a fibered 3-manifolds over the circle with fibered surface $S$ can be seen as the mapping torus of some $\phi\in Mod_\pm(S)$ called its \textit{monodromy}.\par
\textbf{4.1.1. Complement of the figure-eight knot.} In an earlier paper \cite{BR15}, Bridson and Reid have proved that the complement of the figure-eight knot $\mathbb{S}^3\backslash\mathcal{K}$ is profinitely rigid. We denote $F_r$ to be the free group with $r$ generators. Note that the $\mathbb{S}^3\backslash\mathcal{K}$ is a one-punctured torus bundle over the circle with monodromy
\begin{displaymath}
    \psi=\begin{pmatrix}
        2 & 1 \\ 1 & 1
    \end{pmatrix}\in GL_2(\mathbb{Z}). 
\end{displaymath}
And its fundamental group is $F^2\rtimes_\phi\mathbb{Z}$, an extension of the free group of rank 2 by $\mathbb{Z}$. Therefore it is necessary to provide a series of propositions on free-by-cyclic groups. To be precise, they proved:
\begin{proposition}[\cite{BR15}, Theorem B]
    Let $M$ be a compact connected 3-manifold and let $\Gamma=F_r\rtimes_\phi\mathbb{Z}$ be a free-by-cyclic group for some $r\in\mathbb{N}_+$ and $\phi\in Out(F_r)$. If $b_1(M)=1$ and $\widehat{\Gamma}\cong\widehat{\pi_1(M)}$, then $M$ has non-empty boundary, fibers over the circle with compact fiber, and $\pi_1(M)\cong F_r\rtimes\mathbb{Z}$ for some $\psi\in Out(F_r)$.
\end{proposition}
which implies that a 3-manifold whose fundamental group is profinitely isomorphic to a free-by-cyclic group (in particular, a 3-manifold which is profinitely isomorphic to $\mathbb{S}^3\backslash\mathcal{K}$) is itself fibered. Note that any knot complement $X$ has $b_1(X)=1$. 
\begin{proposition}[\cite{BR15}, Lemma 3.1]
    Let $\Gamma_1=N_1\rtimes\mathbb{Z}$ and $\Gamma_2=N_2\rtimes\mathbb{Z}$, with $N_1$ and $N_2$ finitely generated. If $b_1(\Gamma_1)=1$ and $\widehat{\Gamma_1}=\widehat{\Gamma_2}$, then $\widehat{N_1}\cong\widehat{N_2}$.
\end{proposition}
Thus, a 3-manifold profinitely isomorphic to $\mathbb{S}^3\backslash\mathcal{K}$ is fibered with fiber surface the one-punctured torus as well. Finally, the next proposition completes the argument. 
\begin{proposition}[\cite{BR15}, Proposition 3.2]
    Let $F_2$ be the free group of rank 2. For any $\phi\in Out(F_2)=GL_2(\mathbb{Z})$, define $\Gamma_\phi$ to be the free by cyclic group $F^2\rtimes_\phi\mathbb{Z}$. Then for any $\phi_1,\phi_2\in GL_2(\mathbb{Z})$ such that $\widehat{\Gamma_{\phi_1}}=\widehat{\Gamma_{\phi_2}}$, if $\phi_1$ is hyperbolic, then $\phi_2$ is hyperbolic and has the same eigenvalues as $\phi_1$. Equivalently, $det\phi_1=det\phi_2$ and $tr\phi_1=tr\phi_2$.
\end{proposition}
Recall that an elements $\phi$ in $GL_2(\mathbb{Z})$ is \textit{hyperbolic} if it has two different eigenvalues or, equivalently, if the one-punctured torus bundle with monodromy $\phi$ yields a hyperbolic geometry. Since up to conjugate $\begin{pmatrix} 2 & 1 \\ 1 & 1 \end{pmatrix}$ is the unique matrix with trace $3$ and determinant $1$ (easy to check), combination of the three propositions above implies that $\mathbb{S}^3\backslash\mathcal{K}$ is profinitely rigid.\par
\textbf{4.1.2. One-punctured torus bundles over the circle.} The techniques used above is still inspiring and worth a special notion, although soon after that the two authors of \cite{BR15} together with Wilton proved a stronger result in \cite{BRW17} that any one-punctured torus bundle over the circle is profinitely rigid, using deeper concepts about congruence omnipotence in an outer automorphism group. In fact, \cite{BRW17} proves that congruence omnipotence enables one to deduce profinite rigidity results for mapping torus groups $G\rtimes\mathbb{Z}$. We begin with some basic settings from Section 2 of \cite{BRW17}.
\begin{definition}
    \rm{Let $G$ be a finitely generated group and $H\subseteq Out(G)$ be a subgroup. A finite quotient $H\rightarrow Q$ is called a \textit{$G$-congruence quotient} if it factors through $\pi:H\rightarrow P\subseteq Out(G/K)$ where $K$ is a characteristic subgroup of finite index in $G$ and $\pi$ is the restriction of the natural map $Out(G)\rightarrow Out(G/K)$. We say that $Out(G)$ \textit{has the congruence subgroup property} if every finite quotient of $Out(G)$ is a $G$-congruent quotient. More generally, we say that a subgroup $H\subseteq Out(G)$ \textit{has the $G$-congruence subgroup property} if every finite quotient of $H$ is a $G$-congruence quotient. }
\end{definition}
Note that there may be distinct groups $G_1$ and $G_2$ with isomorphic outer automorphism groups $Out(G_1)\cong Out(G_2)$ such that every finite quotient is congruence with respect to $G_1$ but not with respect to $G_2$. For instance, this phenomenon occurs with $Out(F_2)\cong GL_2(\mathbb{Z})\cong Out(\mathbb{Z}^2)$, which has the congruence subgroup property with respect to $F_2$ but not $\mathbb{Z}^2$, Thus "$Out(G)$ has the congruence subgroup property" is a statement about $G$ but not the abstract group $Out(G)$. 
\begin{definition}
    \rm{Let $\Gamma$ be a group. Elements $\gamma_1,\gamma_2\in\Gamma$ of infinite order are said to be \textit{independent} if no non-zero power of $\gamma_1$ is conjugate to a non-zero power of $\gamma_2$ in $\Gamma$. An $m$-tuple $(\gamma_1,\dots,\gamma_m)$ of elements is said to be \textit{independent} if $\gamma_i$ and $\gamma_j$ are independent whenever $1\leqslant i<j\leqslant m$. The group $\Gamma$ is said to be \textit{omnipotent} if, for every independent $m$-tuple $(\gamma_1,\dots,\gamma_m)$ of elements in $\Gamma$, there exists a positive integer $\kappa$ such that, for every $m$-tuple of positive integers $(e_1,\dots,e_m)$ there is a homomorphism to a finite group 
    \begin{displaymath}
        q:\Gamma\rightarrow Q
    \end{displaymath}
    such that $o(q(\gamma_i))=\kappa e_i$, for $i=1,\dots,m$, where $o(g)$ denote the order of a group element $g$. For a subgroup $H$ of $\Gamma$, if we wish to emphasize that an $m$-tuple of elements is independent in $H$, we will say that the tuple is \textit{$H$-independent}.}
\end{definition}
We can finally define the concept of congruent omnipotence. 
\begin{definition}
    \rm{Let $G$ be a finitely generated group, let $H$ be a subgroup of $Out(G)$ and let $\mathcal{S}$ be a subset of $H$. We say that $\mathcal{S}$ is \textit{$(G,H)$-congruence omnipotent} if, for every $m$ and every $H$-independent $m$-tuple $(\phi_1,\dots,\phi_m)$ of elements of $\mathcal{S}$, there is a constant $\kappa$ such that, for any $m$-tuple of positive integers $(n_1,\dots,n_m)$, there is a $G$-congruence quotient $q:H\rightarrow Q$ such that $o(q(\phi_i))=\kappa n_i$ for all $i$.} 
\end{definition}
\cite{BRW17} then proved a general result concerning semi-products over $\mathbb{Z}$.
\begin{theorem}[\cite{BRW17}, Theorem 2.4]
    Let $G$ be a finitely generated group, let $H\subseteq Out(G)$ be a subgroup and let $\mathcal{S}$ be a $(G,H)$-congruence omnipotent subset. Let $\phi_1,\phi_2\in\mathcal{S}$, let $\Gamma_i=G\rtimes_{\phi_i}\mathbb{Z}$ and suppose that $b_1(\Gamma_i)=1$ for $i=1,2$. If $\widehat{\Gamma_1}=\widehat{\Gamma_2}$ then there is an integer $n$ such that ${\phi_1}^n$ is conjugate in $H$ to ${\phi_2}^n$.
\end{theorem}
\cite{Asa01} proved that $GL_2(\mathbb{Z})$ has the congruence subgroup property, while \cite{BER11} proved that virtually free groups, in particular$GL_2(\mathbb{Z})$, are omnipotent. Combination of the two results provides: 
\begin{proposition}[\cite{BRW17}, Proposition 2.7]
    The set of elements of infinite order in $Out(F_2)$ is $F_2$-congruence omnipotent.
\end{proposition}
Then assume there exists a 3-manifold $N$ such that $\widehat{\pi_1M}\cong\widehat{\pi_1N}$. Proposition 4.1 and 4.2 tell us that $N$ is also an one-punctured torus bundle with monodromy, say, $\psi\in GL_2(\mathbb{Z})$. When $\phi$ is not hyperbolic, they use a series of classical algebraic results (ref. \cite{BRW17}, Appendix A) to prove that $\phi$ is conjugate to $\psi$ in $GL_2(\mathbb{Z})$. Thus $\pi_1M\cong\Gamma_\phi\cong\Gamma_\psi\cong\pi_1N$ and we are done. On the other hand, if $\phi$ (respectively, $M$) is hyperbolic, Proposition 4.3 shows that $\psi$ (respectively, $N$) is also hyperbolic and hence $b_1(M)=b_1(N)=1$ by a criterion in \cite{BR15}, Corollary 3.6. Then by \cite{BRW17}, Proposition 2.7 mentioned above, they invoked Theorem 4.7 to deduce that there is an $n$ such that $\phi^n=\psi^{\pm n}$. An algebraic argument is then used to verify that $\phi$ and $\psi$ are conjugate (see the proof of \cite{BRW17}, Theorem A for more details). 
The techniques and procedure introduced above are powerful but still restricting, since they depend heavily on the first betti-number one condition and the explicit algebraic structure of $Out(F_2)=GL_2(\mathbb{Z})$. \par
\textbf{4.1.3. Four-punctured sphere bundles over the circle.} Much more recently, T. Cheetham-West \cite{Che24} proved the profinite rigidity for four-punctured sphere bundles. One of the most remarkable point is that he invoked techniques from \cite{Liu23} (which is our main reference for Section 4.2) to remove the first betti-number one restriction of the methods of Bridson, Reid and Wilton. The main technical result is as below. 
\begin{proposition}[\cite{Che24}, Theorem 1.2]
    Let $M$ ans $N$ be finite-volume hyperbolic manifolds with $\widehat{\pi_1M}\cong\widehat{\pi_1(N)}$. The isomorphism $H^1(M,\mathbb{Z})\rightarrow H^1(M,\mathbb{Z})$ induced by an isomorphism $\Phi:\widehat{\pi_1M}\rightarrow\widehat{\pi_1N}$ sends fibered classes to fibered classes where the corresponding fiber surfaces have the same topological type.  
\end{proposition}
whose proof depends heavily on \cite{Liu23}. Note that profinite isomorphism between groups induces an isomorphism between their (co-)homology, see, for example, \cite{Rei18}. \par
Then just as the procedure of the one-punctured torus bundle case (see the discussion before Proposition 4.8), Cheetham-West proved that $Mod(S^4\backslash4pts)$ is omnipotent and has the congruence subgroup property using its explicit presentation. Next, he proved that any 3-manifold profinitely isomorphic to a four-punctured sphere bundle over the circle is also a four-punctured sphere bundle over a circle, using Proposition 4.9 above and some results in \cite{Liu23}. Also note that \cite{WltZ17} (see Remark 3.4 and Theorem 3.5) distinguishes hyperbolic geometry from others. Finally, a similar justification as the one-punctured torus case (see the discussion after Proposition 4.8) using congruence omnipotence shows that any profinitely isomorphic pair of four-punctured sphere bundle over the circle have conjugate monodromies, hence are homeomorphic. To be precise, we record the statement below.
\begin{theorem}[\cite{Che24}, Lemma 5.5]
    For hyperbolic 3-manifolds $M$ and $N$ that fibered over $S^1$ with $S^3\backslash4pts$ fiber and monodromies $\phi$ and $\psi$ respectively, let $\Phi:\widehat{\pi_1M}\rightarrow\widehat{\pi_1N}$ be an isomorphism that identifies the closures of the fiber subgroups of $M$ and $N$, The monodromies $\phi$ and $\psi$ are conjugate elements of $Mod(S^3\backslash4pts)$, and therefore $M$ and $N$ are homeomorphic. 
\end{theorem}
\textbf{4.1.4. Seeking for examples with higher genus. }Unfortunately, there has not been any profinitely rigid examples of fibered bundles over the circle with a genus $\geqslant
2$ fibered surface. Even the seemingly simplest $2T^2$ case. One of the main problems is that techniques in knot theory cannot be used anymore, while the complexity of presentations of the mapping class groups also matters, see \cite{FarM12}, Chapter 5 for more information. We will return to this topic later.
\subsection{Complete results in almost rigidity}
For any finitely generated residually finite group $\Gamma$, we say that it is \textit{profinitely almost rigid} if there is a finite set of groups $\mathcal{G}$ such that for any group $\Gamma'$ with $\widehat{\Gamma}\cong\widehat{\Gamma'}$, there exists a $\Gamma''\in\mathcal{G}$, such that $\Gamma'\cong\Gamma''$. In other word, there are only a finite number of possibilities for a group profinitely isomorphic to $\Gamma$. Similarly, we say that a 3-manifolds $M$ is profinitely almost rigid if there is a finite set of 3-manifolds $\mathcal{M}$ such that for any 3-manifold $M'$ with $\widehat{\pi_1M}\cong\widehat{\pi_1M'}$, there exists an $M''\in\mathcal{M}$, such that $\pi_1M'\cong\pi_1M''$.\par
After Liu \cite{Liu23} proved that finite-volume hyperbolic 3-manifolds are profinitely almost rigid, it is known that any connected compact orientable 3-manifold is profinitely almost rigid, as concluded in \cite{Xu25}, Table 1 and Table 2. Below we survey the results and techniques invoked in the process.\par
\textbf{4.2.1. Profinite almost rigidity of finite-volume hyperbolic manifolds. }Now we start to survey \cite{Liu23} for a spiritual understanding. The first step is to prove that closed fibered hyperbolic manifolds are profinitely almost rigid. Suppose $M$ is a closed fibered hyperbolic manifold with fibered surface $S$ ($S$ closed with $\chi(S)<0$). Then there is a $\phi\in Mod(S)$ such that $M$ is the mapping torus of $S$, that is, 
\begin{displaymath}
     M=M_f=S\times\mathbb{R}/(x,r+1)\sim(f(x),r), 
\end{displaymath}
where $f$ is a representative of $\phi$. It is known that hyperbolicity of $M$ implies that $\phi$ is pseudo-Anosov, that is, there are two $\phi$-invariant singular measured foliations $(\mathcal{F}_s,\mu_s)$ and $(\mathcal{F}_u,\mu_u)$ respectively on $M$ called the \textit{stable} and \textit{unstable foliations} of $\phi$ respectively, which intersect transversely, such that for some representative $f$ of $\phi$ we have $f(\mathcal{F}_s,\mu_s)=(\mathcal{F}_s,\lambda^{-1}\mu_s)$ and $f(\mathcal{F}_u,\mu_u)=(\mathcal{F}_u,\lambda \mu_u)$ for some $\lambda>1$ determined uniquely by $\phi$ called its \textit{stretch factor}. For a detailed introduction on the classification of mapping class groups and pseudo-Anosov theory, see \cite{FarM12}, Part 3. We also provide some concepts about cohomology. For any mapping torus $M_f$ as above, the \textit{suspension flow} on $M_f$ refers to the continuous family of homeomorphisms $\theta_t:M_f\rightarrow M_f$, parametrized by $t\in\mathbb{R}$ and determined by $(x,r)\mapsto(x,r+t)$. The \textit{distinguished cohomology class} of $M_f$ refers to the cohomology class $\phi_f\in H^1(M_f,\mathbb{Z})$ represented by the distinguished projection $M_f\rightarrow S^1$, where $S^1$ is identified with $\mathbb{R}/\mathbb{Z}$ and the projection is determined by $(x,r)\mapsto r$ (see \cite{Liu23}, Remark 2.4).\par
Recall in Proposition 4.1 and Proposition 4.2 we have already seen that any 3-manifold isomorphic to a fibered manifold is also fibered and has the same fiber surface. Therefor, in order to prove the profinite rigidity of closed fibered 3-manifolds, we need only research the case that $\widehat{\pi_1M}\cong\widehat{\pi_1N}$, $M=M_f$, $N=M_g$ for two pseudo-Anosov elements $f,g$ in $Mod(S)$. The main result in \cite{Liu23} is:
\begin{theorem}[\cite{Liu23}, Corollary 8.2]
    With the same assumption above, the stretch factor of $f$ and $g$ is the same, namely, $\lambda(f)=\lambda(g)$.
\end{theorem}
Then the profinite almost rigidity of closed fibered hyperbolic manifolds follows immediately by the classical result below about pseudo-Anosov mapping classes.
\begin{theorem}[\cite{Liu23}, Theorem 9.2; \cite{FarM12}, Theorem 14.9]
    For any integers $g,p\geqslant0$ with $2-2g-p<0$, and for any real constant $C>1$, there are at most finitely many pseudo-Anosov automorphisms $f:S\rightarrow S$ of an orientable surface $S$ with genus $g$ and $p$ punctures, up to topological equivalence, such that the stretch factor $\lambda(f)$ is at most $C$.
\end{theorem}
Before we sketch the verification of Theorem 4.11, which is quite skillful and needs some new concepts, we first clarify how the result in the fibered case leads to the general case. The key observation is that if a finite cover of a 3-manifold $M$ is profinitely almost rigid, so is $M$ itself (\cite{Liu23}, Lemma 9.3). But virtually compact specialization of hyperbolic manifolds implies that they are virtually fibered (see Appendix A for an introduction and also \cite{AscFW15}, Flowchart 4). Thus the case of closed hyperbolic manifolds follows. Finally, the argument in the proof of \cite{Liu23}, Lemma 9.5 says that it extends without much difficulty to the bounded case.\par
Now we only have to concentrate on the argument toward Theorem 4.11. We will need some setting about periodic points and their indices (ref. \cite{Liu23}, Section 2 and Section 8). Suppose that $M$ is an orientable connected closed fibered 3–manifold that admits a hyperbolic metric. Suppose that $\phi\in H^1(M,\mathbb{Z})$ is a primitive fibered class. Recall that we say $\phi$ is a \textit{fibered cohomology class} if there exists some fibration of $M$ with fiber surface $S$ and pseudo-Anosov monodromy $f$ (that is, $M=M_f$), such that $\phi$ is the distinguished cohomology class $\phi_f$. In this setting, periodic orbits of $f$ correspond to periodic trajectories of the suspension flow, and they have well-defined periodic indices.\par
To be precise, for any natural number $m\in\mathbb{N}$, an \textit{$m$-periodic point} of $f$ refers to a fixed point $p\in S$ of $f^m$, and it gives rise to an $m$-periodic trajectory of the suspension flow on $M_f$, namely, the loop $\mathbb{R}/m\mathbb{Z}\rightarrow M_f$ determined by the map $\mathbb{R}\rightarrow S\times\mathbb{R}:r\mapsto(p,r)$. The free-homotopy class of an $m$-periodic trajectory depends only on the $f$-iteration orbit of an $m$-periodic point. We denote by $Orb_m(f)$ the set of $m$-periodic orbits of $f$, and for any $m$-periodic orbit $\mathbf{O}\in Orb_m(f)$, denote by $l_m(f;\mathbf{O})$ the free-homotopy class of the $m$-periodic trajectory determined by $\mathbf{O}$. Since the homotopy set $[\mathbb{R}/m\mathbb{Z},M_f]$ is naturally identified with the set of conjugacy classes of $\pi_1(M_f)$, denoted as $Orb(\pi_1(M_f))$, we write
\begin{displaymath}
    l_m(f;\mathbf{O})\in Orb(\pi_1(M_f)),
\end{displaymath}
referring to it as an essential \textit{$m$-periodic trajectory class}.\par
For any $m$–periodic point $p\in S$ of $f$, the $m$–periodic index of $f$ at $p$ refers to the fixed point index of $f^m$ at $p$. Explicitly, when $p$ is a $k$–prong singularity of the stable (or unstable) invariant foliation for some $k\geqslant 3$, the $m$–periodic index of $f$ at $p$ equals $1-k$ if $f^m$ preserves every prong of the foliation at $p$, or it equals $1$ otherwise; when $p$ is a regular point, the $m$–periodic index of $f$ at $p$ equals $-1$ or $1$, according to a similar rule as if $k=2$. Note that the $m$–periodic index of $f$ depends only on the $f$–iteration orbit of $m$–periodic points. For any $m$–periodic orbit $\mathbf{O}\in Orb_m(f)$ , we define the $m$–periodic index of $f$ at $\mathbf{O}$ using any point $p\in\mathbf{O}$, and denote it as
\begin{displaymath}
    ind_m(f,\mathbf{O})\in\mathbb{Z}\backslash\{0\}.
\end{displaymath}\par
Now for any $m\in \mathbb{N}$ and $i\in\mathbb{Z}\backslash\{0\}$, we denote the number of index-$i$ $m$-periodic orbits of $f$ as
\begin{displaymath}
    \nu_m(M,\phi;i)=\nu_m(f;i)=\#\{\mathbf{O}\in Orb_m(f):ind_m(f;\mathbf{O})=i\},
\end{displaymath}
and denote the number of (essential) $m$-periodic orbits of $f$ as
\begin{displaymath}
    N_m(M,\phi)=N_m(f)\sum_{i\in\mathbb{Z}\backslash\{0\}}\nu_m(f;i)=\#Orb_m(f). 
\end{displaymath}
The number $N_m(f)$ is called the \textit{$m$-th orbit Nielsen number} of $f$. It is known that the orbit Nielsen numbers growth exponentially fast as $m$ tends to $\infty$, and we have the following alternative representation of the stretch factor
\begin{displaymath}
    \lambda(f)=\limsup_{m\rightarrow\infty}N_m(f)^{1/m}, 
\end{displaymath}
see \cite{Jia96}, Section 2, Example 2.\par
Now let $M_A,M_B$ be a pair of profinitely isomorphic closed fibered hyperbolic 3-manifolds whose monodromies are $f,g$ respectively. Denote $\Psi$ to be the isomorphism $\widehat{\pi_1M_A}\rightarrow\widehat{\pi_1M_B}$. Liu proved in \cite{Liu23}, Theorem 8.1 that there are a pair of properly chosen (see below) fibered classes $\phi_B\in H^1(M_B,\mathbb{Z})$ and $\psi_A=\Psi^*_{1/\mu}(\phi_B)\in H^1(M_A,\mathbb{Z})$ such that
\begin{displaymath}
    \nu_m(M_A,\psi_A;i)=\nu_m(M_B,\phi_B;i)
\end{displaymath}
for any $m\in\mathbb{N}$ and $i\in\mathbb{Z}\backslash\{0\}$. Hence, for any $m\in\mathbb{N}$, 
\begin{displaymath}
    N_m(M_A,\psi_A)=N_m(M_B,\phi_B), 
\end{displaymath}
and the representation of the stretch factor above implies $\lambda(f)=\lambda(g)$ and we have proved Theorem 4.11. \par
The final task is to give a properly choice of the pair of fibered classes $(\psi_A,\phi_B)$. Denote $H_A=H_1(M_A,\mathbb{Z})$ and $H_B=H_1(M_B,\mathbb{Z})$. Then $\Psi$ induces an isomorphism on the level of abelization $\Psi_*:\widehat{H_A}\rightarrow\widehat{H_B}$. Identify $\widehat{H_A}=H_A\otimes_\mathbb{Z}\widehat{\mathbb{Z}}$ and $\widehat{H_B}=H_B\otimes_\mathbb{Z}\widehat{\mathbb{Z}}$. We define the \textit{matrix coefficient module} for $\Phi$ to be the smallest $\mathbb{Z}$-submodule $L$ of $\widehat{\mathbb{Z}}$, such that $\Psi_*(H_A)$ lies in the submodule $H_B\otimes_\mathbb{Z}L$ of $\widehat{H_B}=H_B\otimes_\mathbb{Z}\widehat{\mathbb{Z}}$ (here we see $H_A$ as a submodule of $\widehat{H_A}$ since $H_A$ is residually finite, see Section 2). We denote the matrix coefficient module as $MC(\Psi_*)$. The next proposition says that $MC(\Psi_*)$ is of dimension $1$ in $\widehat{\mathbb{Z}}$.
\begin{proposition}[\cite{Liu23}, Theorem 6.1]
    There exists some unit $\mu\in\widehat{\mathbb{Z}}^\times$, such that $MC(\Psi_*)$ is the $\mathbb{Z}$-submodule $\mu\mathbb{Z}$ of $\widehat{\mathbb{Z}}$.
\end{proposition}
Note that the scalar multiplication by $1/\mu$ on $\widehat{\mathbb{Z}}$ defines a homomorphism $1/\mu\in Hom_\mathbb{Z}(MC(\Psi_*),\mathbb{R})$ with image in $\mathbb{Z}$. Therefore, we obtain a homomorphism $\Psi_{*,1/\mu}:H_A\otimes_\mathbb{Z}\mathbb{R}\rightarrow H_B\otimes_\mathbb{Z}\mathbb{R}$ by the unique R–linear extension of the composition homomorphism
\begin{displaymath}
    H_A\stackrel{\Psi_*|_{H_A}}{\longrightarrow}H_B\otimes_\mathbb{Z}MC(\Psi_*)\stackrel{1\otimes1/\mu}{\longrightarrow}H_B\otimes_\mathbb{Z}\mathbb{R}.
\end{displaymath}
called the $1/\mu$-specialization of $\Psi_*$ and denoted as $\Psi_*^{1/\mu}$. A duality of homology and cohomology gives a homomorphism on cohomology
\begin{displaymath}
    \Psi^*_{1/\mu}:H^1(M_B;\mathbb{R})\rightarrow H^1(M_A,\mathbb{R}).
\end{displaymath}
\cite{Liu23}, Theorem 5.1 (see also Theorem 5.16 of this paper) then says that if $\Psi^*_{1/\mu}$ is nondegenerate, then it witnesses a dimension-preserving bijective correspondence between the Thurston-norm cones for $M_A$ and those for $M_B$. Namely, under the linear isomorphism $\Psi^*_{1/\mu}$, every Thurston-norm cone for $M_B$ projects onto a distinct and unique Thurston-norm cone for $M_A$. Now the classical result concerning the correspondence between Thurston-norm cones and fibered cones (see \cite{Thu86} for an original reading and Subsection 5.2 for an introduction) shows the existence of the pair $(\psi_A,\phi_B)$ with $\psi_A=\Psi^*_{1/\mu}(\phi_B)$. 
\begin{remark}
    \rm{In the process of arguments in \cite{Liu23}, there is also an application of the combinative torsion, such as some kind of Lefchetz function and twisted Reidemeister torsions. We do not include this part in this paper and would like to refer the reader to \cite{Liu23}, Section 7 and relative sections.}
\end{remark}
\textbf{4.2.2. Profinite almost rigidity of 3-manifolds.} \cite{Xu25} concludes that any connected compact orientable 3-manifold with empty or toroidal boundary is profinitely rigid. We follow Table 1 and Table 2 of \cite{Xu25} to give references for non-hyperbolic cases.\par
Profinite almost rigidity of $Sol$-manifolds is proved by Grunewald, Pickel and Segal in \cite{GPS80}. Profinite almost rigidity of closed Seifert fibered manifolds is proved by Wilkes in \cite{Wilk17}. Profinite almost rigidity of closed graph manifolds is proved by Wilkes in \cite{Wilk18}. Profinite almost rigidity of mixed and bounded non-hyperbolic manifolds is proved by Xu in \cite{Xu25}. And finally profinite almost rigidity of finite-volume hyperbolic manifolds is proved by Liu in \cite{Liu23}. 
\subsection{The work of Agol, Cheetham-West and Minsky}
However, when we really come to the research of profinite rigidity of hyperbolic manifolds, there exists much more inevitable difficulties as one may imagine. For example, we do not know whether the profinite rigidity pf a finite index regular covering of a 3-manifold $M$ can lead to the profinite rigidity of $M$ (not like the case of profinite almost rigidity, see \cite{Liu23}, Lemma 9.3). Nevertheless, much more recently Agol, Cheetham-West and Minsky proved in \cite{ACM24} a first reduction to the problem. Namely, they proved the next theorem. Note that a group is the fundamental group of a finite-volume hyperbolic manifold if and only if it is a lattice in $PSL_2(\mathbb{C})$, the isometry group of $\mathbb{H}^3$.
\begin{theorem}[\cite{ACM24}, Corollary 4.7]
   The following are equivalent:\par
   (1) The profinite rigidity of all non-arithmetic lattices in $PSL_2(\mathbb{C})$.\par
   (2) The profinite rigidity of all fibered non-arithmetic lattices in $PSL_2(\mathbb{C})$.\par
   (3) The profinite rigidity of all special non-arithmetic lattices in $PSL_2(\mathbb{C})$.
\end{theorem}
By $(1)\Leftrightarrow(2)$ of the above theorem, we need only consider the profinite rigidity of fibered hyperbolic manifolds and arithmetic hyperbolic manifolds. There has not been satisfactory discussion on each part, though.\par
For the complete argument the reader is referred to the original paper for a detailed reading since it is brief enough and has a more group-theoretic taste compared to other results of this note. It is based closely on the Malnormal Special Quotient Theorem in \cite{Wis12a} and virtual compact specialization of hyperbolic manifolds (which also has a close connection to the former). See Appendix A for a brief introduction of virtual specialization of irreducible 3-manifolds. Results on omnipotence in \cite{She23} are also invoked, see \cite{ACM24}, Theorem 2.4. The main technical result in \cite{ACM24} is the following.
\begin{proposition}[\cite{ACM24}, Theorem 4.1 and 4.2]
    For $\Gamma_1$ a torsion-free non-arithmetic lattice in $PSL_2(\mathbb{C})$, and $\Delta<\Gamma_1$ a finite-index subgroup, there is a lattice $\Gamma_2<\Delta<\Gamma$ such that $Aut^+(\Gamma_2)=\Gamma_1$. 
\end{proposition}
Here $Aut^+(G)$ denotes the orientation-preserving automorphisms of the lattice $G$.\par
We sketch the proof that $(2)\Rightarrow(1)$ in Theorem 4.15. For a non-arithmetic lattice $\Gamma<PSL_2(\mathbb{C})$ we can choose $\Delta<\Gamma$ such that $\Delta$ is the fundamental group of a hyperbolic 3-manifold that fibers over the circle, using virtual specialization. Then it will follow that the finite-index subgroup $\Theta<\Delta$ furnished by the propositions will be a fibered lattice as well, satisfying $Aut^+(\Theta)=\Gamma$. By \cite{BR22}, Theorem 4.4, once we assume that $\Theta$ is profinitely rigid, $\Gamma$ will be profinitely rigid as well. The proof that $(3)\Rightarrow(2)$ is similar. 
\section{Summary of methods}
In this section we summarize some techniques using in the research process above. Most of them have the same favor that we need to extend classical results in the original setting of 3-manifold topology. A great example is the "profinite hyperbolization theorem" (Theorem 3.5) and the classical Hyperbolization Theorem (Theorem 3.9). The former states that hyperbolicity of a 3-manifold is equivalent to absence of $\widehat{\mathbb{Z}}^2$ subgroup of its profinite fundamental group, while the latter says the same thing in the original setting, that hyperbolicity of a 3-manifold is (almost) equivalent to absence of essential embedded torus, that is, $\mathbb{Z}^2$ subgroup of its fundamental group. 
\subsection{The classical Bass-Serre theory and its profinite version}
We begin with a brief review of the classical Bass-Serre theory and its profinite version. We have claimed its importance in the profinite detection of JSJ decomposition in Section 3. Our main reference for this part is \cite{Xu25}, Section 3 and \cite{WltZ19}.\par
An \textit{oriented graph} $\Gamma$ is a disjoint union $E(\Gamma)\sqcup V(\Gamma)$ of sets which are called its \textit{vertex set} and \textit{edge set} respectively, with two \textit{relation maps} $d_0,d_1:E(\Gamma)\rightarrow V(\Gamma)$ that assign the \textit{initial vertex} and \textit{terminal vertex} of each edge. A \textit{morphism} $\alpha:\Gamma\rightarrow\Delta$ of oriented graphs is a continuous map with $\alpha d_i=d_i\alpha$ for $i=0,1$.\par
A \textit{graph of groups}, denoted by $(\mathcal{G},\Gamma)$, consists of a connected oriented graph $\Gamma$, a group $\mathcal{G}_x$ associated to each $x\in\Gamma=V(\Gamma)\sqcup E(\Gamma)$, and an injective homomorphism $\phi^i_e:\mathcal{G}_e\rightarrow\mathcal{G}_{d_i(e)}$ for each $e\in E(\Gamma)$ and $i=0,1$.\par
A geometric interpretation of graph of groups turns out to be a graph of topological spaces. Indeed, let $X_v$ (respectively, $X_e$) be CW-complexes so that $\pi_1X_v=\mathcal{G}_v$ (respectively, $\pi_1X_e=\mathcal{G}_e$), and let $\phi^e_i:X_e\rightarrow X_{d_i(e)}$ be the maps (up to homotopy) so that $\phi^e_i$ induces 
\begin{displaymath}
    \varphi^e_i:\pi_1X_e=\mathcal{G}_e\rightarrow\pi_1X_{d_i(e)}=\mathcal{G}_{d_i(e)}.
\end{displaymath}
We can construct a CW-complex according to these information. 
\begin{displaymath}
    X_\Gamma=(\bigsqcup_{e\in E(\Gamma)}X_e\times[0,1])\bigcup_{\phi^e_i:X_e\times\{i\}\rightarrow X_{d_i}(e)}(\bigsqcup_{v\in V(\Gamma)}X_v)
\end{displaymath}\par
It turns out that $\pi_1(X_\Gamma)\cong\Pi_1^{abs}(\mathcal{G},\Gamma)$, see \cite{SW79}. 
\begin{definition}
    \rm{Let $(\mathcal{G},\Gamma)$ be a finite graph of groups and $S$ be a maximal subtree of $\Gamma$, the \textit{(abstract) fundamental group} with respect to $S$ is defined as
    \begin{displaymath}
        \Pi^{abs}_1(\mathcal{G},\Gamma,S)=\frac{(*_{v\in V(\Gamma)}\mathcal{G}_v)*F_{E(\Gamma)}}{the\ normal\ closure\ of\ \{t_e:e\in E(S)\}\sqcup\{\phi_1^e(g)^{-1}t_e\phi_0^e(g)t_e^{-1}:e\in E(\Gamma)\ and\ g\in\mathcal{G}_e\}}, 
    \end{displaymath}
    where $F_{E(\Gamma)}$ is the free group over $E(\Gamma)$ with generators $t_e$, $e\in E(\Gamma)$.}
\end{definition}
It follows from \cite{Ser80}, Proposition 20 that up to isomorphism, $\Pi^{abs}_1(\mathcal{G},\Gamma,S)$ is independent with the choice of the maximal subtree $S$. Thus, we may omit the notation $S$ and simply denote the fundamental group as $\Pi^{abs}_1(\mathcal{G},\Gamma)$ when there is no need for a precise maximal subtree.\par
Turning to the profinite version, we begin with some definitions.
\begin{definition}
    \rm{A \textit{profinite graph} $\Gamma$ is a graph such that: \par
    (1) $\Gamma$ is a profinite space (that is, an inverse limit of finite discrete spaces);\par
    (2) $V(\Gamma)$ is closed;\par
    (3) the maps $d_0$ and $d_1$ are continuous.}
\end{definition}
By \cite{RibZ10}, Proposition 2.1.4 every profinite graph $\Gamma$ is an inverse limit of finite quotient graphs of $\Gamma$.\par
For a profinite space $X$ that is the inverse limit $X_j$, $[[\widehat{\mathbb{Z}}X]]$ is defined to be the inverse limit of $[\widehat{\mathbb{Z}}X]$ which is the free $\widehat{\mathbb{Z}}$-module with basis $X_j$. For a pointed profinite space $(X,*)$ that is the inverse limit of pointed finite discrete spaces $(X_j,*)$, $[[\widehat{\mathbb{Z}}(X,*)]]$ is the inverse limit of $[\widehat{\mathbb{Z}}(X_j,*)]$, where $[\widehat{\mathbb{Z}}(X_j,*)]$ is the free $\widehat{\mathbb{Z}}$-module with basis $X_j\backslash\{*\}$ (see \cite{RibZ10}, Chapter 5.2).\par
For a graph $\Gamma$ define the pointed space $(E^*(\Gamma),*)$ as $\Gamma/V(\Gamma)$ with the image of $V(\Gamma)$ as a distinguished point $*$, and denote the image of $e\in E(\Gamma)$ by $\overline{e}$.
\begin{definition}
    \rm{A \textit{profinite tree} $\Gamma$ is a profinite graph such that the sequence
    \begin{displaymath}
        0\rightarrow[[\widehat{\mathbb{Z}}(E^*(\Gamma),*)]]\stackrel{\delta}{\rightarrow}[[\widehat{\mathbb{Z}}V(\Gamma)]]\stackrel{\epsilon}{\rightarrow}\widehat{\mathbb{Z}}\rightarrow 0
    \end{displaymath}
    is exact, where $\delta(\overline{e})=d_1(e)-d_0(e)$ for every $e\in E(\Gamma)$ and $\epsilon(v)=1$ for every $v\in V(\Gamma)$.}
\end{definition}
We call any graph of groups a \textit{graph of profinite groups} if its vertex and edge groups are all profinite groups. The profinite fundamental group $\Pi_1(\mathcal{G},\Gamma)$ of a finite graph of finitely generated profinite groups $(\mathcal{G},\Gamma)$ can be defined as the profinite completion of its abstract (usual) fundamental group (see Definition 5.1), denoted as $\Pi_1^{abs}(\mathcal{G},\Gamma)$. The profinite fundamental group $\Pi_1(\mathcal{G},\Gamma)$ has the following presentation:
\begin{displaymath}
    \Pi_1(\mathcal{G},\Gamma)=\langle\mathcal{G}_v,t_e\ |\ rel(\mathcal{G}_v),d_1(g)=t_ed_0(g)t_e^{-1},g\in\mathcal{G}_e,t_e=1\ for\ e\in T\rangle;
\end{displaymath}
where $rel(\mathcal{G}_v)$ is the set of relaters of $\mathcal{G}_v$ and $T$ is any maximal subtree of $\Gamma$.\par
It is not difficult to see that the vertex groups of $(\mathcal{G},\Gamma)$ embed in $\Pi_1^{abs}(\mathcal{G},\Gamma)$ but not always embed in $\Pi_1(\mathcal{G},\Gamma)$. IF they do embed, $(\mathcal{G},\Gamma)$ is called \textit{injective}. If $(\mathcal{G},\Gamma)$ is not injective, the edge and vertex groups can be replaced by their images in $\Pi_1(\mathcal{G},\Gamma)$, and after this replacement $(\mathcal{G},\Gamma)$ becomes injective. See \cite{Rib17}, Section 6.4 for a detailed introduction.\par
Now we turn to the tree structure and for convenience we state the concept of a \textit{profinite Bass-Serre tree} directly. 
\begin{definition}
    \rm{Let $(\mathcal{G},\Gamma)$ be a finite graph of profinite groups, and denote $\Pi=\Pi_1(\mathcal{G},\Gamma)$. For each $v\in V(\Gamma)$, let $\Pi_v$ denote the image of $\mathcal{G}_v$ in $\Pi$ through the map $\mathcal{G}_v\rightarrow(\sqcup_{v\in V(\Gamma)}\mathcal{G}_v)\sqcup\mathcal{F}_{E(\Gamma)}\rightarrow\Pi$ where $\mathcal{F}_{E(\Gamma)}$ is the free profinite group over the finite space $E(\Gamma)$ with generators denoted by $t_e\ (e\in E(\Gamma))$ (see the definition of $\Pi_1(\mathcal{G},\Gamma)$ above), and for each $e\in E(\Gamma)$, let $\Pi_e$ denote the image of $\mathcal{G}_e$ in $\Pi$ through the relation map $\phi_1^e:\mathcal{G}_e\rightarrow\mathcal{G}_{d_1(e)}$. Then, for each $x\in V(\Gamma)\cup E(\Gamma)$, $\Pi_x$ is a closed subgroup in $\Pi$, so the coset space $\Pi/\Pi_v$ is a profinite space. Let
    \begin{displaymath}
        V(T)=\bigsqcup_{v\in V(\Gamma)}\Pi/\Pi_v,\ E(\Gamma)=\bigsqcup_{e\in E(\Gamma)}\Pi/\Pi_e,\ and\ T=V(T)\sqcup E(T)
    \end{displaymath}
    be equipped with disjoint-union topology. The relation map is given by
    \begin{displaymath}
        d_0(g\Pi_e)=gt_e\Pi_{d_0(e)},\ d_1(g\Pi_e)=g\Pi_{d_1(e)},\ \forall e\in E(\Gamma),g\in\Pi.
    \end{displaymath}
    Then $T$ is called the \textit{profinite Bass-Serre tree} associated to $(\mathcal{G},\Gamma)$}
\end{definition}
Note that $T$ is a disjoint of finitely many profinite spaces. Hence $T$ is a profinite space. It is easy to verify that the relation maps $d_0,d_1$ is continuous. Moreover, \cite{Rib17}, Corollary 6.3.6 shows that $T$ is indeed a profinite tree (see Definition 5.3).\par
Finally we turn to group action on graphs. We also describe it on the profinite settings directly. By definition a profinite group $G$ \textit{acts} on a profinite graph $\Gamma$ if we have a continuous action of $G$ on the profinite space $\Gamma$ that commutes with the maps $d_0$ and $d_1$. We include the next lemma.\par
\begin{lemma}
    Suppose that a profinite group $G$ acts on a profinite tree $T$ and does not fix any vertex. Then there exists an open normal subgroup $U$ of $G$ that is not generated by its vertex stabilizers. 
\end{lemma}
Now as an application of the discussion above, we will restate Theorem 3.3 in an stronger and exact form (see \cite{WltZ19}, Theorem 4.3) and sketch the proof there. Recall that for any irreducible 3-manifold $M$, the union of its JSJ tori of its JSJ decomposition forms a submanifold, denoted as $\mathcal{T}$.  Then $\mathcal{T}$ induces a graph-of-spaces decomposition of $M$ and a graph-of groups decomposition of $\pi_1M$, then a graph-of-profinite groups decomposition of $\widehat{\pi_1M}$. Denote the Bass-Serre tree of the last decomposition to be $\widehat{T}_M$ (note that we add the hat in order to distinguish it from the ordinary Bass-Serre tree of $\pi_1M$, although there is no need to introduce it in our discussion). Wilton and Zalesskii then proved the much stronger theorem below. 
\begin{theorem}[\cite{WltZ19}, Theorem 4.3]
    If $M,M'$ are closed orientable, irreducible 3-manifolds and 
    \begin{displaymath}
        f:\widehat{\pi_1M}\stackrel{\cong}{\rightarrow}\widehat{\pi_1M'}
    \end{displaymath}
    is an isomorphism, then there is an $f$-equivariant isomorphism
    \begin{displaymath}
        \phi:\widehat{T}_M\rightarrow\widehat{T}_M'
    \end{displaymath}
    of the corresponding profinite Bass-Serre trees. In particular, the underlying graphs of the JSJ decompositions of $M$ and $M'$ are isomorphic, as are the profinite completions of the fundamental groups of the corresponding pieces. 
\end{theorem}
The idea to prove Theorem 5.6 is to consider the action of the vertex and edge groups of $M$ (that is, fundamental groups of JSJ pieces and JSJ tori, respectively) on $\widehat{T}_{M'}$. First, they proved that for any vertex space $N$ of $M$ (that is, $N$ is a JSJ piece of $M$), $\widehat{\pi_1N}$ must act with a point on $\widehat{T}_{M'}$. The proof of this statement occupies Lemma 4.4-4.6 of \cite{WltZ19} by treating hyperbolic pieces and Seifert fibered pieces separately. Then, they proved some similar (but more complex) results corresponding to edge groups of $M$ (that is, $\widehat{Z}^2$ subgroups of $\widehat{\pi_1M}$), which occupies Lemma 4.7-4.9 of \cite{WltZ19}.\par
Now we can construct a map $\phi:\widehat{T}_M\rightarrow\widehat{T}_{M'}$. To be precise, they proved:
\begin{proposition}[\cite{WltZ19}, Lemma 4.10]
    Consider closed, orientable, irreducible 3-manifolds $M,M'$ and let $f:\widehat{\pi_1M}\rightarrow\widehat{\pi_1M'}$ be an isomorphism. Then there exists an $f$-equivariant morphism of graphs
    \begin{displaymath}
        \phi:\widehat{T}_M\rightarrow\widehat{T}_{M'}.
    \end{displaymath}
    Note that, here, we only claim that $\phi$ is a map of abstract, non-oriented graphs. This map may in principle send edges to either edges or vertices. 
\end{proposition}
Finally, note that by symmetry we have the morphism in the opposite direction
\begin{displaymath}
    \psi:\widehat{T}_{M'}\rightarrow\widehat{T}_M,
\end{displaymath}
which is $f^{-1}$-equivariant. Equivariance implies that
\begin{displaymath}
    g\psi\circ\phi(x)=\psi\circ\phi(gx)
\end{displaymath}
for all $g\in\widehat{\pi_1M}$ and $x\in\widehat{T}_M$, whence the stabilizer of $x$ is contained in the stabilizer of $\psi\circ\phi(x)$. Since vertex-stabilizers stabilize unique vertices, it follows that $\psi\circ\phi$ is equal to the identity on vertices, and hence on the whole of $\widehat{T}_M$. In particular, $\phi$ and $\psi$ induce isomorphisms of the finite quotient graphs $\widehat{\pi_1M}/\widehat{T}_M$ and $\widehat{\pi_1M'}/\widehat{T}_{M'}$. We may therefore choose consistent orientations on these graphs, which lift to equivariant orientations on the profinite trees $\widehat{T}_M$ and $\widehat{T}_{M'}$, which are respected by $\phi$ and $\psi$. This also implies continuity of $\phi$ and $\psi$. The proof of Theorem 4.10 is finally accomplished.
\subsection{Thurston norm and Thurston norm cones}
In this subsection we introduce a series of concepts and results used significantly in \cite{Liu23}, that is, the Thurston norm defined on the homology and cohomology of a 3-manifold. For convenience, we only state the results for a closed 3-manifold. It is not difficult to extend them to the bounded case. We refer the reader to Thurston's original paper \cite{Thu86} and also \cite{Liu23}, Section 2 and Section 6 for a more detailed introduction.\par
Thurston proved the following theorem for the construction of the semi-norm.
\begin{theorem}[\cite{Thu86}, Theorem 1]
    Let $M$ be a closed oriented 3-manifold with possibly empty boundary. There is a canonical continuous function $\lVert\cdot\rVert_{Th}$ defined on the second homology groups of $M$, $H_2(M;\mathbb{R})$, which is convex and linear on rays through the origin.\par
    If, furthermore, every embedded surface representing a non-zero element of ($H_2(M;\mathbb{Z})$ has negative Euler characteristic, then $\lVert\cdot\rVert_{Th}$ is a norm.\par
    In general, $\lVert\cdot\rVert_{Th}$ is a semi-norm vanishing on precisely the subspace spanned by embedded surfaces of non-negative Euler characteristic. 
\end{theorem}
Before we give the explicit construction, note that using Poincar\'e duality, we carry $\lVert\cdot\rVert_{Th}$ over to a function on $H^1(M;\mathbb{R})$. Furthermore, When $\lVert\cdot\rVert_{Th}$ is a norm on $H_2(M)$, we have also a dual norm $\lVert\cdot\rVert^*_{Th}$ on the dual vector spaces $H_1(M;\mathbb{R})=H^2(M;\mathbb{R})$, defined by the formula
\begin{displaymath}
    \lVert\alpha\rVert^*_{Th}=\sup_{\lVert\beta\rVert_{Th}\leqslant1}\{\alpha\cdot\beta\}.
\end{displaymath}\par
Now we begin to verify Theorem 5.8. The method is similar to the one how we solve Cauchy's functional equation $f(x+y)=f(x)+f(y)$ for continuous $f$. We define $\lVert\cdot\rVert_{Th}$ first on the integral lattice of $H_2(M;\mathbb{R})$ by the formula
\begin{displaymath}
    \lVert a\rVert_{Th}=\inf\{\chi_\_(S)\ |\ S\ is\ an\ embedded\ surface\ representing\ a\}.
\end{displaymath}\par
We invoke the following elementary lemma.
\begin{lemma}[\cite{Thu86}, Lemma 1]
    In an closed oriented 3-manifold $M$, every element $a\in H_2(M;\mathbb{Z})$ is represented by an embedded oriented surface $S$. If $a$ is divisible by some $k\in\mathbb{N}_+$, then $S$ is a union of $k$ components, each representing $a/k$.
\end{lemma}
Then it follows that for every integer $k$, $\lVert ka\rVert_{Th}\geqslant k\lVert a\rVert_{Th}$. On the other hand, if a surface $S$ represents $a$, then $k$ parallel copies represent $ka$, implying that $\lVert ka\rVert_{Th}=k\lVert a\rVert_{Th}$. The linearity is done.\par
The triangle inequality $\lVert a+b\rVert_{Th}\leqslant \lVert a\rVert_{Th}+\lVert b\rVert_{Th}$ is then proved by writing the representing surfaces of $a$ and $b$ and using an Euler characteristic argument. We have proved that $\lVert\cdot\rVert_{Th}$ is a semi-norm on the integral lattices.\par
Next we extend $\lVert\cdot\rVert_{Th}$ to the rational points in $H_2(M;\mathbb{R})$ by the condition that $\lVert\cdot\rVert_{Th}$ is linear on each ray through the origin. By the inequality above, $\lVert\cdot\rVert_{Th}$ is a convex function. It follows easily that we can extend $\lVert\cdot\rVert_{Th}$ to all of $H_2(M;\mathbb{R})$ in a unique way so that it is continuous. The extended function $\lVert\cdot\rVert_{Th}$ is still convex, and linear on rays.\par
The statement concerning the zero space of $\lVert\cdot\rVert_{Th}$ is then easily verified and we are done.\par
Next we turn to some objects corresponding to the Thurston norm. We begin with the unit ball. Note that any norm $\lVert\cdot\rVert$ is determined by its unit ball
\begin{displaymath}
    \mathcal{B}=\{x\ |\ \lVert x\rVert\leqslant1\}.
\end{displaymath}
Our norms $\lVert\cdot\rVert_{Th}$ and $\lVert\cdot\rVert_{Th}^*$ are not like the more familiar norms derived from inner products, since, as we shall see, their unit balls, denoted as $\mathcal{B}_{Th}(M)$ and $\mathcal{B}^*_{Th}(M)$ respectively, are (possibly noncompact) codimension-$0$ convex polyhedra rather than ellipsoids. 
\begin{theorem}[\cite{Thu86}, Theorem 2]
    When $\lVert\cdot\rVert_{Th}$ is a norm, the unit ball $\mathcal{B}^*_{Th}(M)$ is a polyhedron whose vertices are lattice points, $\pm\beta_1,\cdots,\pm\beta_k$ and the unit ball $\mathcal{B}_{Th}(M)$ is a polyhedron defined by the linear inequalities with integer coefficients
    \begin{displaymath}
        \mathcal{B}_{Th}(M)=\{a\ |\ \lvert a\cdot\beta_i\rvert\leqslant1\ (1\leqslant i\leqslant k)\ \}
    \end{displaymath}
\end{theorem}
For the proof and some examples of computation, see \cite{Thu86}, Section 2.\par
It is worthwhile to note that in the definition of the Thurston norm, we may use properly immersed subsurfaces instead of properly embedded ones, and the resulting semi-norm does not change. This nontrivial fact is proved by Gabai in \cite{Gab83}, confirming a former conjecture of Thurston in \cite{Thu86}, Section 5. In particular, the Thurston norm for any finite cover is proportional to the covering degree.
\begin{theorem}[\cite{Gab83}, Corollary 6.13]
    Let $M$ be any orientable connected compact 3-manifold, and $M'\rightarrow M$ be a finite cover. Then the following formula holds for any cohomology class $\phi\in H^1(M;\mathbb{R})$ and its pullback $\phi'\in H^1(M';\mathbb{R})$:
    \begin{displaymath}
        \lVert\phi'\rVert_{Th}=[M':M]\cdot\lVert\phi\rVert_{Th}.
    \end{displaymath}
\end{theorem}
Now we turn to the most important case, closed fibered 3-manifolds, which produces a series of concepts about faces of the norm unit ball and cones above them. We begin with some basic setting. Suppose that $M$ is a closed 3-manifold that fibers over a circle. Then each possible fiber is a incompressible embedded surface. And any incompressible surface in the homology class of a fiber is homotopic to a fiber. This can be seen by passing to the infinite cyclic covering induced by projection to $S^1$. Any surface homologous to a fiber lifts to this cover, which is homeomorphic to (a fiber)$\times\mathbb{R}$. The projection of a lift to a fiber is a degree $1$ map which is injective on $\pi_1$, hence is homotopic to a homeomorphism. In fact, \cite{Thu86}, Theorem 4 shows that any incompressible surface in the homology class of a fiber is isotopic to a fiber.\par
Now consider the space $\tau$ of tangent planes to the fibers. If we pick an orientation for a fiber, and use this to orient $\tau$, it makes sense to speak of the Euler class $\chi(\tau)$ (the first obstruction to a section of $\tau$) as an element of $H^2(M;\mathbb{Z})$. We have the following: 
\begin{proposition}[\cite{Thu86}, Theorem 3]
    Let $M$ be a compact oriented 3-manifold which fibers over $S^1$ with fiber a surface of negative Euler characteristic. The Euler class of the tangent space to any fibration of $M$ over $S^1$ is a vertex of $\mathcal{B}^*_{Th}(M)$ and in particular
    \begin{displaymath}
        \lVert\chi(\tau)\rVert_{Th}^*=1.
    \end{displaymath}
    The ray determined by the homology class of any fiber passes through the interior of a top-dimensional face of $\partial \mathcal{B}_{Th}(M)$, and the formula
    \begin{displaymath}
        \lVert a\rVert_{Th}=|\chi(\tau)\cdot a|
    \end{displaymath}
    holds in some neighborhood of this ray. 
\end{proposition}
By the proposition above the following corollary is easily verified.
\begin{corollary}
    If $\lVert\cdot\rVert_{Th}$ is a norm in $H_2(M;\mathbb{R})$ and if $H_2(M;\mathbb{R})$ has rank $\geqslant2$, then $M$ possesses at least one incompressible surface which is not the fiber of a fibration.
\end{corollary}
Since $\mathcal{B}_{Th}(M)$ and $\mathcal{B}_{Th}^*(M)$ are all polyhedra, it is natural to consider their faces and cones above them. Namely, we have the next:
\begin{proposition}[\cite{Thu86}, Theorem 5]
    The set $F$ of cohomology classes in $H^1(M;\mathbb{R})$ representable by non-singular closed $1$-form is some union of the cones on open faces of $\mathcal{B}_{Th}(M)$, minus the origin. The set of elements in $H^1(M;\mathbb{Z})$ whose Lefchetz dual is represented by a fiber of a fibration consists of the set of all lattice points in $F$.
\end{proposition}
Recall that a cohomology class $\phi\in H^1(M;\mathbb{Z})$ is said to be \textit{fibered} if $\phi$ is represented by the free homotopy class of a continuous map $M\rightarrow S^1$ which is a fiber bundle projection onto the standard circle. The fiber is connected precisely when $\phi$ is primitive. The above proposition by Thurston shows that every fibered class $\phi$ lies in the interior of a unique codimension-$0$ Thurston norm cone, denoted by $\mathcal{C}_{Th}(M,\phi)$, and moreover, every cohomology class $\psi\in H^1(M;\mathbb{Z})$ in the interior of $\mathcal{C}_{Th}(M,\phi)$ is a fibered class. Any $\mathcal{C}_{Th}(M,\phi)$ arising this way is called a \textit{fibered cone} for $M$.\par
Next we consider the cones on $H_1(M;\mathbb{R})$ dual to fibered cones in $H^1(M;\mathbb{R})$. These homology cones are intimately related to flow structures on $M$ transverse to the surface fibers. We describe the picture below focusing on the case of fibered closed hyperbolic 3-manifolds. See \cite{FLP12}, Exposition 14 for a more detailed introduction.\par
Suppose that $M$ is an orientable closed 3-manifold which admits a hyperbolic metric, $\phi\in H^1(M;\mathbb{Z})$ is a primitive fibered class. Identify $M$ with the mapping torus of a pseudo-Anosov automorphism such that $\phi$ is identified with its distinguished cohomology class (see Subsection 4.2). D. Fried gives a characterization of the linear dual in $H_1(M;\mathbb{R})$ of the fibered cone $\mathcal{C}_{Th}(M,\phi)$, which we denote as
\begin{displaymath}
    \mathcal{C}^{Fr}(M,\phi)=\{x\in H_1(M;\mathbb{R}):\psi(x)\geqslant0\ for\ all\ \psi\in\mathcal{C}_{Th}(M,\phi)\}.
\end{displaymath}\par
Considering all the homology classes represented by the forward periodic trajectories of the suspension flow, Fried shows that the radical rays in $H_1(M;\mathbb{R})$ that pass through those homology classes form a dense subset of $\mathcal{C}^{Fr}(M\phi)$. Therefore, $\mathcal{C}^{Fr}(M,\phi)$ is the smallest convex polyhedral cone in $H_1(M;\mathbb{R})$ that contains those homology classes, and it has codimension zero in $H_1(M;\mathbb{R})$.\par
Then we turn to the projective world of homology (as vector spaces). For any orientable closed 3-manifold $M$, we denote by $\mathbf{P}(H_1(M;\mathbb{R}))$ the projectivization of $H_1(M;\mathbb{R})$. The points of $\mathbf{P}(H_1(M;\mathbb{R}))$ are considered to be the real 1-dimensional linear subspaces of $H_1(M;\mathbb{R})$. Tor any codimension-0 cone $\mathcal{C}$ in $H^1(M;\mathbb{R})$ with respect to the Thurston norm of $M$, we introduce a subset of $\mathbf{P}(H_1(M;\mathbb{R}))$ as follows:
\begin{displaymath}
    \mathcal{D}(M,\mathcal{C})=\{l\in\mathbf{P}(H_1(M;\mathbb{R})):
    l\cap \ker(\phi)=\{0\}\ for\ all\ \phi\in int(\mathcal{C})\}.
\end{displaymath}\par
In general, $\mathcal{D}(M,\mathcal{C})$ is projectively isomorphic to a polytope. Its codimension in $\mathbf{P}(H_1(M;\mathbb{R}))$ equals the kernel dimension of the Thurston norm. We refer to $\mathcal{D}(M,\mathcal{C})$ as the \textit{projective dual polytope} with respect to $\mathcal{C}$.\par
The main result proved in \cite{Liu23} concerning the projective dual polytope is the following, which describe the behavior of it under finite covering.\par
\begin{proposition}[\cite{Liu23}, Lemma 6.7]
    Let $M$ be an orientable connected closed 3–manifold which admits a hyperbolic metric. If $M$ has a fibered cone $\mathcal{C}$, then there exists a finite cover $M'$ of $M$, such that $\mathcal{D}(M',\mathcal{C}')$ has a vertex which projects into the interior of $\mathcal{D}(M,\mathcal{C})$.
\end{proposition}
We record the proof of the proposition from \cite{Liu23} as an exposition of the techniques. A similar treatment also occurs in \cite{Liu20}, Problem 1.5 and Outline of the Solution. Fix a primitive cohomology class $\phi\in H^1(M,\mathbb{Z})$ in $\mathcal{C}$, we identify $(M,\phi)$ with the mapping torus of a pseudo-Anosov automorphism $f:S\rightarrow S$ and its distinguished cohomology class, namely, $(M_f,\phi_f)$. Since $M_f$ is closed and hyperbolic, $\pi_1(M_f)$ is nonelementary word-hyperbolic and virtually compact special (see \cite{Ago13}, Theorem 9.3 and for more about the virtually compact specialization see subsection 5.3). By Wise's special quotient theorem (see \cite{Wis12b}, Theorem 12.7), there exists a group quotient $\pi_1(M_f)\rightarrow G$, such that $G$ is again nonelementary word-hyperbolic and virtually compact special. In particular, $G$ has some finite-index subgroup $G'$ with $b_1(G')>0$ (see \cite{Ago13}, Corollary 1.2).\par
Now we can take $M'$ to be the finite cover of $M_f$ such that $\pi_1(M')$ is the preimage of $G'$ with respect to $\pi_1(M_f)\rightarrow G$. The pullback of $\phi_f$ is a (not necessarily primitive) fibered class $\phi'\in H^1(M';\mathbb{Z})$, and it lies in a unique fibered cone $\mathcal{C}'$, which contains the image of $\mathcal{C}=\mathcal{C}_f$. The suspension flow on $M_f$ lifts to a flow on $M'$. If one identifies $M'$ with the mapping torus of a pseudo-Anosov automorphism $f':S'\rightarrow S'$, where $S'$ if a preimage component of $S$, then the lifted flow is just the suspension flow with velocity scaled down by the divisibility of $\phi'$.\par
For any vertex $v'$ of $\mathcal{D}(M',\mathcal{C}')$, there exists a periodic trajectory $\gamma'$ of the lifted flow, such that $[\gamma']\in H_1(M';\mathbb{R})$ lies in the linear 1–subspace $v'$, by Fried's dual characterization of $\mathcal{C}'$. The periodic trajectory $\gamma'$ covers a unique periodic trajectory $\gamma$ in $M_f$. The free-homotopy class of $\gamma$ can be represented as the conjugacy class of some element $g$ in $\pi_1(M_f)$. If the vertex $v'$ projects into a codimension–1 closed face of $\mathcal{D}(M_f,\mathcal{C}_f)$, the image of $g$ in $G$ has finite order, by the above construction and Lemma 6.6 of \cite{Liu23}. Then, in view of the commutative diagram
\begin{displaymath}
    \begin{tikzcd}
        H_1(M';\mathbb{R}) \arrow[d] \arrow[r] & H_1(G';\mathbb{R}) \arrow[d]\\
        H_1(M_f;\mathbb{R}) \arrow[r] & H_1(G;\mathbb{R})
\end{tikzcd}
\end{displaymath}
the 1-subspace $v'$ of $H_1(M';\mathbb{R})$ must have trivial image in $H_1(G;\mathbb{R})$. Note that the homomorphisms in the above diagram are all surjective.\par
By Fried's characterization, $\mathcal{D}(M',\mathcal{C}')$ has codimension $0$ in $\mathbf{P}(H_1(M';\mathbb{R}))$, and $\mathcal{D}(M_f,\mathcal{C}_f)$ has codimension $0$ in $\mathbf{P}(H_1(M_f;\mathbb{R}))$. Since $b_1(G)>0$, there must be some vertex $w'$ of $\mathcal{D}(M',\mathcal{C}')$, and as an $1$-subspace of $H_1(M';\mathbb{R})$, $w'$ has nonvanishing image in $H_1(G;\mathbb{R})$. The argument in the last paragraph shows that $w'$ must project into the interior of $\mathcal{D}(M_f,\mathcal{C}_f)$. Therefore, $M'$ is a finite cover as asserted. We have proved Proposition 5.15.\par
Finally, We state a result by Liu answering how profinite isomorphic of a 3-manifold pair determines the correspondence of their Thurston norm cones. We have referred to it also in the sketch of proof of Theorem 4.11
\begin{theorem}[\cite{Liu23}, Theorem 5.1]
    Let $M_A,M_B$ be irreducible 3-manifolds with empty or tori boundary which admit complete metrics of non-positive curvature in the interior and $\Psi:\widehat{\pi_A}\rightarrow\widehat{\pi_B}$ be a profinite isomorphism of their fundamental groups. If $\Psi^*_\epsilon:H^1(M_B;\mathbb{R})\rightarrow H^1(M_A;\mathbb{R})$ is non-degenerate, then $\Psi^*_\epsilon$ witnesses a dimension-preserving bijective correspondence between the Thurston-norm cones for $M_A$ and those for $M_B$. Namely, under the linear isomorphism $\Psi^*_\epsilon$, every Thurston-norm cone for $M_B$ projects onto a distinct and unique Thurston-norm cone for $M_A$.
\end{theorem}
Note that the homomorphism $\Psi^*_\epsilon$ is defined in the discussion after Proposition 4.13. (there we define it for $\epsilon=1/\mu$).\par
Now, as an application of the discussion above, we sketch the proof of the closed case of Proposition 4.13. Recall that we let $M_A,M_B$ be closed hyperbolic 3-manifolds of positive first Betti number and $\pi_A,\pi_B$ are their fundamental groups, respectively. $\Psi$ is an isomorphism $\widehat{\pi_A}\stackrel{\cong}{\rightarrow}\widehat{\pi_B}$. We want to prove that there exists some unit $\mu\in\widehat{\mathbb{Z}}^{\times}$ such that $MC(\Psi_*)$ is the $\mathbb{Z}$-submodule $\mu\mathbb{Z}$ of $\widehat{\mathbb{Z}}$.\par
First, by Theorem 5.11 and the virtually fibering of hyperbolic manifolds, we may assume that $M_A,M_B$ fiber over the circle. We make the following constructions. Fix a fibered cone $\mathcal{C}_A$ in $H^1(M_A;\mathbb{R})$. Obtain a finite cover $M'_A$ of $M_A$ by Proposition 5.15. Thus, $\mathcal{D}(M_A',\mathcal{C}_A')$ has some vertex $v_A'$ that projects into the interior of $\mathcal{D}(M_A,\mathcal{C}_A)$, where $\mathcal{C}_A'$ refers to the unique fibered cone for $M_A'$ that contains the image of $\mathcal{C}_A$, with respect to the covering-induced homomorphism $H^1(M_A;\mathbb{R})\rightarrow H^1(M_A';\mathbb{R})$. Denote by $q_A$ the image of $v_A'$ in the interior of $\mathcal{D}(M_A,\mathcal{C}_A)$. Obtain the finite cover $M_B'\rightarrow M_B$ which makes a $\Psi$-corresponding pair with $M_A'\rightarrow M_A$, that is, such that $\pi_A'=\pi_A\cap\Psi^{-1}(clos(\pi_B'))$ where $clos(\pi_B')$ denotes the closure of $\pi_B'$ in $\widehat{\pi_B}$. Liu then argued a pair of fibered cones $(\mathcal{C}_B,\mathcal{C}_B')$ "corresponding" to $(\mathcal{C}_A,\mathcal{C}_A')$, namely, the following diagram
\begin{displaymath}
    \begin{tikzcd}
        \mathcal{D}(M_A',\mathcal{C}_A') \arrow[d,"\mathbf{P}(cov_*)"] \arrow[r,"\mathbf{P}((\Psi')^{\epsilon'}_*)"] & \mathcal{D}(M_B',\mathcal{C}_B') \arrow[d,"\mathbf{P}(cov_*)"]\\
        \mathcal{D}(M_A,\mathcal{C}_A) \arrow[r,"\mathbf{P}(\Psi^\epsilon_*)"] & \mathcal{D}(M_B,\mathcal{C}_B)
\end{tikzcd}
\end{displaymath}
commutes for a certain open set of $\epsilon$ in $\mathbb{R}$ and a corresponding open set of $\epsilon'$, where $\mathbf{P}(cov_*)$ is the projectivization of the homology homomorphism induced by coverings and $\mathbf{P}(\Psi^\epsilon_*),\mathbf{P}((\Psi')^{\epsilon'}_*)$ are projectivization of the corresponding homomorphisms defined in the discussion after Proposition 4.13. (there we define them for $\epsilon=1/\mu$). Then Liu defined $v_B'=\mathbf{P}((\Psi')^{\epsilon'}_*)(v_A')$ and $v_B=\mathbf{P}(\Psi^\epsilon_*)(v_A)$ and showed that
the diagram
\begin{displaymath}
    \begin{tikzcd}
        v_A' \arrow[d] \arrow[r] & v_B' \arrow[d]\\
        q_A \arrow[r] & q_B
    \end{tikzcd}
\end{displaymath}
of specified elements is also commutative. The conclusion then follows from a technical lemma, \cite{Liu23}, Lemma 6.4.
\appendix
\section{Virtual specialization}
After Perelman proved the Geometrization Theorem in 2002, one of the most important achievement in 3-manifold topology may be the series of results concerning virtual (compact) specialization of 3-manifolds. These results obtain their importance by the fact that virtually special 3-manifolds are virtually fibered (see \cite{AscFW15}, Flowchart 4). Virtually fibering is used significantly in the research of profinite properties of 3-manifolds, for example, Subsection 4.2 and the proof of Proposition 5.15. Therefore, we would like to briefly review this important topic in the last subsection of this paper. See \cite{AscFW15}, Section 4 for a more detailed introduction.\par
Recall that a 3-manifold is called \textit{(compact) special} if its fundamental group is isomorphic to the fundamental group of a non-positively curved (compact) special cube complex. For the exact definition see also \cite{HW07}.\par
For irreducible 3-manifolds according to their JSJ decompositions we classify them into three classes, hyperbolic manifolds, mixed manifolds and graph manifolds. We begin with non-hyperbolic cases as well.\par
For mixed manifolds, Przytycki and Wise \cite{PW12} proved that they are virtually special. For graph manifolds, Liu \cite{Liu13} proved that a non-positively curved graph manifold is virtually special. And it is an early result of Leeb \cite{Leb95} that graph manifolds with non-empty boundary are non-positively curved. Then graph manifolds with non-empty boundary are virtually special. This result is also proved independently by Przytycki and Wise \cite{PW14}. As for closed graph manifolds, Buyalo and Svetlov \cite{BS05} discussed the criteria to judge when they are non-positively curved. Finally, Hagen and Przytycki \cite{HP13} provided an example to show that non-positively curved graph manifolds may not be virtually compact special.\par
Virtual compact specialization of closed hyperbolic manifolds is proved by Agol \cite{Ago13}, depending on earlier work of Kahn-Markovic \cite{KM12} and Sageev-Wise \cite{SaW12}, while virtual compact specialization of cusped hyperbolic manifolds is proved by Wise \cite{Wis12a}.\par


\begin{thebibliography}{99}
\bibitem[Ago13]{Ago13}I. Agol, \textit{The virtual Haken conjecture}, with an appendix by I. Agol, D. Groves, J. Manning, Documenta Math. 18 (2013), 1045–1087.
\bibitem[ACM24]{ACM24}I. Ago, T. Cheetham-West and Y. Minsky, \textit{Simply transitive geodesics and omnipotence of lattices in $PSL(2,\mathbb{C})$}, arXiv:2409.08418v1, 2024.
\bibitem[Asa01]{Asa01}M. Asada, \textit{The faithfulness of the monodromy representations associated with certain families of algebraic curves}, J. Pure Appl. Algebra, 159(2-3):123–147, 2001.
\bibitem[AscFW15]{AscFW15}M. Aschenbrenner, S. Friedl, and H. Wilton, \textit{3-Manifold Groups}, EMS Series of Lectures in Mathematics, 2015.
\bibitem[BER11]{BER11}K. Bux, M. V. Ershov, and A. S. Rapinchuk, \textit{The congruence subgroup property for $Aut(F_2)$: a group-theoretic proof of Asada's theorem}, Groups Geom. Dyn., 5(2):327–353, 2011.
\bibitem[BR15]{BR15}M. R. Bridson, A. W. Reid, \textit{Profinite rigidity, fibering, and the figure-eight knot}, arXiv:1505.07886v1, 2015.
\bibitem[BR22]{BR22}M. R. Bridson and A. W. Reid, \textit{Profinite rigidity, Kleinian groups and the cofinite Hopf property}, Michigan Math. J. 72 (2022), 25–49.
\bibitem[BRW17]{BRW17}M. R. Bridson, A. W. Reid and H. Wilton, \textit{Profinite rigidity and surface bundles over the circle}, Bull. London Math. Soc. 49 (2017) 831841.
\bibitem[BS05]{BS05}S. Buyalo, P. Svetlov, \textit{Topological and geometric properties of graph-manifolds}, St. Petersburg Math. J. 16 (2005), no. 2, 297–340.
\bibitem[Che24]{Che24}T. Cheetham-West, Profinite rigidity and hyperbolic four-punctured sphere bundles over the circle, arXiv:2308.00266v3, 2024.
\bibitem[FarM12]{FarM12}B. Farb and D. Margalit, \textit{A Primer on Mapping Class Groups}, Princeton Mathematical Series, 49. Princeton University Press, Princeton, NJ, 2012.
\bibitem[FLP12]{FLP12}A. Fathi, F. Laudenbach, V. Poénaru, \textit{Thurston’s Work on Surfaces}, translated from the 1979 French original by D. M. Kim and D. Margalit. Princeton University Press, Princeton, NJ, 2012.
\bibitem[Fun13]{Fun13}L. Funar, \textit{Torus bundles not distinguished by TQFT invariants}, Geometry \& Topology 17 (2013), no. 4, 2289–2344.
\bibitem[Gab83]{Gab83}D. Gabai, \textit{Foliations and the topology of 3–manifolds}, J. Differential Geom. 18 (1983), no. 3, 445–503.
\bibitem[GPS80]{GPS80}F. J. Grunewald, P. F. Pickel, and D. Segal, \textit{Polycyclic groups with isomorphic finite quotients}, Annals of Mathematics 111 (1980), no. 1, 155-195.
\bibitem[Hem87]{Hem87}J. Hempel, \textit{Residual finiteness for 3-manifolds}, Combinatorial Group Theory and Topology, pp. 379–396, Annals of Mathematics Studies, vol. 111, Princeton University Press, Princeton, NJ, 1987.
\bibitem[Hem14]{Hem14}J. Hempel, \textit{Some 3-manifold groups with the same finite quotients}, arXiv:1409.3509.
\bibitem[HP13]{HP13}M. Hagen, P. Przytycki, \textit{Cocompactly cubulated graph manifolds}, Israel J. Math., to appear.
\bibitem[HW07]{HW07}F. Haglund and D.l Wise, \textit{Special cube complexes}, Geom. Funct. Anal. (2007), 1–69.
\bibitem[Jia96]{Jia96}B. Jiang, \textit{Estimation of the number of periodic orbits}, Pacific J. Math. 172 (1996), 151–185.
\bibitem[KM12]{KM12}J. Kahn, V. Markovic, Immersing almost geodesic surfaces in a closed hyperbolic three manifold, Ann. of Math. (2) 175 (2012), 1127–1190.
\bibitem[Leb95]{Leb95}B. Leeb, \textit{3-manifolds with(out) metrics of nonpositive curvature}, Invent. Math. 122 (1995), 277–289.
\bibitem[Liu13]{Liu13}Y. Liu, \textit{Virtual cubulation of nonpositively curved graph manifolds}, J. Topol. 6 (2013), no. 4, 793–822.
\bibitem[Liu20]{Liu20}Y. Liu, \textit{Virtual homological spectral radii for automorphisms of surfaces}, J. Amer. Math. Soc. 33 (2020), no. 4, 1167–1227.
\bibitem[Liu23]{Liu23}Y. Liu, \textit{Finite-volume hyperbolic 3-manifolds are almost determined by their finite quotient groups}, Inventiones mathematicae 231 (2023), no. 2, 741–804.
\bibitem[NS07]{NS07}N. Nikolov and D. Segal, \textit{On finitely generated profinite groups}, I. Strong completeness and uniform bounds. Ann. of Math. (2), 165(1):171–238, 2007.
\bibitem[PW12]{PW12}P. Przytycki, D. Wise, \textit{Mixed 3-manifolds are virtually special}, preprint, 2012.
\bibitem[PW14]{PW14}P. Przytycki, D. Wise, \textit{Graph manifolds with boundary are virtually special}, J. Topology 7 (2014), 419–435.
\bibitem[Rei18]{Rei18}A. W. Reid, \textit{Profinite rigidity}, Proceedings of the International Congress of Mathematicians—Rio de Janeiro 2018. Vol. II. Invited lectures, pp. 1193–1216, World Sci. Publ., Hackensack, NJ, 2018.
\bibitem[Rib17]{Rib17}L. Ribes, \textit{Profinite graphs and groups}, volume 66 of Ergebnisse der Mathematik und ihrer Grenzgebiete. 3. Folge. A Series of Modern Surveys in Mathematics [Results in Mathematics and Related Areas. 3rd Series. A Series of Modern Surveys in Mathematics]. Springer, Cham, 2017.
\bibitem[RibZ10]{RibZ10}L. Ribes and P. Zalesskii, \textit{Profinite groups}, Ergebnisse der Mathematik und ihrer Grenzgebiete (A Series of Modern Surveys in Mathematics), Springer Berlin Heidelberg, 2010.
\bibitem[Ser80]{Ser80}J. Serre, \textit{Trees}, Springer Monographs in Mathematics, Springer Berlin Heidelberg, 1980.
\bibitem[She23]{She23}S. Shepherd, \textit{Imitator homomorphisms for special cube complexes}, Trans. Am. Math. Soc. 376(1) (2023), 599–641.
\bibitem[Ste72]{Ste72}P. F. Stebe, \textit{Conjugacy separability of groups of integer matrices}, Proceedings of the American Mathematical Society 32 (1972), no. 1, 1–7.
\bibitem[SW79]{SW79}P. Scott and T. Wall, \textit{Topological methods in group theory}, London Mathematical Society Lecture Note Series, p. 137–204, Cambridge University Press, 1979.
\bibitem[SaW12]{SaW12}M. Sageev, D. Wise, \textit{Cubing cores for quasiconvex actions}, preprint, 2012.
\bibitem[Thu82]{Thu82}W. P. Thurston, \textit{Three dimensional manifolds, Kleinian groups and hyperbolic geometry}, Bull. Amer. Math. Soc. 6 (1982), 357–379. 
\bibitem[Thu86]{Thu86}W. P. Thurston, \textit{A norm for the homology of 3-manifolds}, Mem. Amer. Math. Soc. 59 (1986), no. 339, pp. 99–130.
\bibitem[Wilk17]{Wilk17}G. Wilkes, \textit{Profinite rigidity for Seifert fiber space}s, Geometriae Dedicata 188 (2017), no. 1, 141-163.
\bibitem[Wilk18]{Wilk18}G. Wilkes, \textit{Profinite rigidity of graph manifolds and JSJ decompositions of 3-manifolds}, Journal of Algebra 502 (2018), 538–587.
\bibitem[Wis12a]{Wis12a}D. Wise, \textit{The structure of groups with a quasi-convex hierarchy}, preprint, 2012.
\bibitem[Wis12b]{Wis12b}D. Wise, \textit{From Riches to Raags: 3-Manifolds, Right-Angled Artin Groups, and Cubical Geometry}, CBMS Regional Conference Series in Mathematics, 117, American Mathematical Society, Providence, RI, 2012.
\bibitem[WltZ10]{WltZ10}H. Wilton and P. Zalesskii, \textit{Profinite properties of graph manifolds}, Geometriae Dedicata, 147(1):29—45, 2010.
\bibitem[WltZ17]{WltZ17}H. Wilton and P. A. Zalesskii, \textit{Distinguishing geometries using finite quotients}, Geometry \& Topology 21 (2017), no. 1, 345-384
\bibitem[WltZ19]{WltZ19}H. Wilton and P. A. Zalesskii, \textit{Profinite detection of 3-manifold decompositions}, Compos. Math. 155 (2019), no. 2, 246–259.
\bibitem[Xu25]{Xu25}X. Xu, \textit{Profinite almost rigidity in 3-manifolds}, arXiv:2410.16002v3, 2025.
\end{thebibliography}
\end{document}